\documentclass[9pt]{amsart}  
\usepackage{amsfonts}
\usepackage{graphicx}
\usepackage{amscd}
\usepackage{amsmath,amsfonts,amssymb,amsthm,mathrsfs,xypic,cite}
\xyoption{curve}
\usepackage{fancybox}
\newcommand{\m}{\Lambda}
\newcommand{\Hom}{\operatorname{Hom}}

\newcommand{\Ker}{\operatorname{Ker}}
\newcommand{\cok}{\operatorname{Coker}}
\newcommand{\Ima}{\operatorname{Im}}

\newcommand{\ha}{\operatorname{\mathcal{A}}}
\newcommand {\hp}{\mathcal P}
\newcommand{\s}{\hfill \blacksquare}
\newtheorem{thm}{Theorem}[section]

\newtheorem{cor}[thm]{Corollary}

\newtheorem{lem}[thm]{Lemma}
\newtheorem{sublem}[thm]{Sub-Lemma}
\newtheorem{exm}[thm]{Example}

\newtheorem{prop}[thm]{Proposition}
\newtheorem{rem}[thm]{Remark}
\newtheorem{defn}[thm]{Definition}
\begin{document}

\title [Monic monomial modules]
{Monic monomial representations ${\rm \bf I:}$\\
Gorenstein-projective modules}
\author [Xiu-Hua Luo, Pu Zhang ] {Xiu-Hua Luo, Pu Zhang$^*$}
\thanks{$^*$The corresponding author.}
\thanks{Supported by the NSFC 11271251, 11431010 and 11401323.}

\begin{abstract} \ For a $k$-algebra $A$, a quiver $Q$, and an ideal $I$ of $kQ$ generated by monomial relations,
let $\Lambda: = A\otimes_k kQ/I$. We introduce the monic
representations of  $(Q, I)$ over $A$. We give properties of the
structural maps of monic representations, and prove that the
category ${\rm mon}(Q, I, A)$ of the monic representations of $(Q,
I)$ over $A$ is a resolving subcategory of ${\rm rep}(Q, I, A)$. We
introduce the condition ${\rm(G)}$. The main result claims that a
$\m$-module is Gorenstein-projective if and only if it is a monic
module satisfying ${\rm(G)}$. As consequences, the monic
$\m$-modules are exactly the projective $\m$-modules if and only if
$A$ is semisimple; and they are exactly the Gorenstein-projective
$\m$-modules if and only if $A$ is selfinjective, and if and only if
${\rm mon}(Q, I, A)$ is Frobenius.

\vskip5pt

\end{abstract} \maketitle
\section {\bf Introduction and reliminaries}

\subsection{} The representation category ${\rm rep}(Q, I, k)$  of a
bounded quiver $(Q, I)$ over field $k$ ([28], [4]), and the morphism
category ([3]) of $k$-algebra $A$, can be reviewed
as the representation category ${\rm rep}(Q, I, A)$ of $(Q, I)$ over
$A$, or equivalently, the module category $\m\mbox{-}{\rm mod}$ with
$\m: = A\otimes_k \ kQ/I$. If $I = 0$, this viewpoint induces the notion of {\it monic}
$\m$-modules ([25]), and we get {\it the
monomorphism category} ${\rm mon}(Q, A)$ consisting of
monic representations of $Q$ over $A$.

\vskip5pt

In fact, there is a long history of studying  ${\rm mon}(Q, A)$,
initiated by G.Birkhoff [7]. When $Q$ is of type $A_n$, it is {\it
the submodule category} ([29]-[31], [26]) or {\it the filtered chain
category} ([1], [33], [34]). In particular, C.M.Ringel and
M.Schmidmeier have established the Auslander-Reiten theory of the
submodule category ([30]); and D.Simson have studied the
representation type ([33]). A reciprocity of monomorphism operator
and perpendicular operator is given ([35]); and by D.Kussin,
H.Lenzing and H.Meltzer ([23], [24]) and X.W.Chen [10], the
monomorphism category is related to the singularity theory.

\vskip5pt However, for  $I\ne 0$, how to define monic
representations of $(Q, I)$ over $A$ is a problem. In this paper we
consider the case that $I$ is generated by monomial relations.
Definition \ref{maindef} is inspired by the property of projective
$kQ/I$-modules. This fits our aim of describing
Gorenstein-projective $\m$-module via monic representations of $(Q,
I)$ over $A$. If $A$ is not semisimple, then the monic
representations of $(Q, I)$ over $A$ are much more than the
projective $\m$-modules. We give properties of the structural maps
$X_p$ of a monic representation $X$ (Theorem \ref{m4}). The category
${\rm mon}(Q, I, A)$ of the monic representations of $(Q, I)$ over
$A$ is proved to be a resolving subcategory of ${\rm rep}(Q, I, A)$
(Theorem \ref{resolving}).

\subsection{} On the other hand, Gorenstein-projective modules (M.Auslander and
M.Bridger [2]; E.E.Enochs and O.M.G.Jenda [13]) enjoy more pleasant
stable properties than projective modules (cf. [8], [6], [32],
[18]). They form a main ingredient in the relative homological
algebra (cf. [14], [5], [17], [11], [20]), and are widely used in
the representation theory and algebraic geometry (cf. [8], [22],
[9], [16], [6], [12], [19]). An important feature is that the
category $\mathcal {GP}(A)$ of Gorenstein-projective modules is
Frobenius, and hence its stable category $\underline{\mathcal
{GP}(A)}$ is triangulated. If $A$ is Gorenstein then the singularity
category of $A$ is triangle equivalent to $\underline{\mathcal
{GP}(A)}$ ([8]), and hence has the Auslander-Reiten triangles
([15]).

\vskip5pt

Thus we need explicitly construct all the Gorenstein-projective
modules of an algebra $\m$. This turns out to be closely related to
the monomorphism category (cf. [35], [19], [25], [36]). However, if
we view an algebra as $\m: = A\otimes_k \ kQ/I$, this relation is only
known at the case of $I = 0$ ([25]).

\vskip5pt

We introduce the condition ${\rm(G)}$. The main result Theorem
\ref{mainthm} claims that a $\m$-module is Gorenstein-projective if
and only if it is a monic representation of $(Q, I)$ over $A$
satisfying the condition ${\rm(G)}$. Note that here $\m$ is not
necessarily Gorenstein. In this way we give an inductive
construction of the Gorenstein-projective $\m$-modules via quivers.
This generalizes the corresponding result in [25] for $I =0$, however,
the proof of this general case is more complicated, and a main tool
is a description of the Gorenstein-projective modules over the
triangular extension of two algebras via a bimodule ([36]). So we
get
\begin{align*}\{\mbox{projective} \ \m\mbox{-module}\} & \ \subseteq
\ \{\mbox{Gorenstein-projective} \ \m\mbox{-module}\} \\ & \ = \
\{\mbox{monic} \ \m\mbox{-module satisfying {\rm(G)}}\} \subseteq \
\{\mbox{monic} \ \m\mbox{-module}\}.\end{align*}As applications of
Theorem \ref{mainthm}, the first inclusion is an equality if and
only if $A$ is semisimple; and the last inclusion is an equality if
and only if $A$ is selfinjective, and if and only if ${\rm mon}(Q,
I, A)$ is Frobenius.

\subsection{} Throughout $A$ is a finite-dimensional algebra over field $k$, $A$-mod the
category of finitely generated left $A$-modules, $Q$ a finite
acyclic quiver, $I$ an ideal of path algebra $kQ$ generated by
monomial relations, and $\m: = A\otimes \ kQ/I$. All tensors
$\otimes$ are over $k$. Let $Q_0$ and $Q_1$ be the set of vertices
and the set of arrows, respectively, $s(\alpha)$ and $e(\alpha)$ the
starting and the ending vertex of arrow $\alpha$, respectively. Path
$p$ is a sequence $\alpha_l\cdots \alpha_1$ of arrows with each
$e(\alpha_i) = s(\alpha_{i+1});$ and it is {\it non-zero} if
$p\notin I$.  Connection of paths is from right to left.  Vertex $v$
is a path of length $0$, denoted by $e_v$, and $P(v)$ is the
indecomposable projective $kQ/I$-module $(kQ/I)e_v$. Let $\mathcal
P$ be the set of paths, $s(p)$, $e(p)$ and $l(p)$ the starting
vertex, the ending vertex and the length of $p\in \mathcal P$,
respectively. Label $Q_0$ as $1, \cdots, n$, such that $j
> i$ if $\alpha:
j\rightarrow i$ is in $Q_1$. So $n$ is a source and $1$ is a sink.
If $p=\alpha q\in\mathcal P$ with $\alpha\in Q_1$, then  $\alpha$ is
{\it the last arrow} of $p$, denoted by ${\rm la}(p): = \alpha$. Let
$\rho: =\{\rho_1, \cdots, \rho_t\}$ is the minimal set of generators
of $I$. For $i, j \in Q_0$, we put
\begin{align*}& \ha(\to i): = \{\alpha\in Q_1 \ | \
e(\alpha) =i\}\\& \hp(\to i): = \{\ p\in\mathcal P \ |  \  e(p) =
i,  \ l(p) \ge 1,  \
 p \notin I\}\\ &\ha(j\to i): = \{ \alpha\in Q_1 \ | \ s(\alpha) = j, \ e(\alpha) =
i\};\\ &\hp(j\to i): = \{\ p\in\mathcal P \ | \ s(p) = j, \ e(p) =
i, \ l(p)\ge 1, \ p \notin I\}.\end{align*} \noindent For a non-zero
path $p$ with $l(p) \ge 1$, we put
$$K_p: = \{q\in \hp(\to s(p)) \ | \ pq \in
I\}. \eqno (1.1)$$

\subsection{} {\it A representation $X$ of
$Q$ over $A$} is a datum $(X_i, \ X_{\alpha}, \ i\in Q_0,  \ \alpha
\in Q_1)$, where each $X_i$ is an $A$-module and each $X_{\alpha}:
X_{s(\alpha)} \rightarrow X_{e(\alpha)}$ is an $A$-map. A morphism $f: X\rightarrow
Y$ is a datum $(f_i, \ i\in Q_0)$, where $f_i: X_i \rightarrow Y_i$
is an $A$-map, such that for each arrow $\alpha: j\rightarrow i$ we
have $$Y_\alpha f_j = f_iX_\alpha. \eqno(1.2)$$ We call $X_i$ and
$f_i$ {\it the $i$-th branch} of $X$ and  {\it the $i$-th branch} of
$f$, respectively. We also write $X$ as $\left(\begin{smallmatrix}
X_1\\
\vdots\\
X_n
\end{smallmatrix}\right)_{(X_\alpha, \ \alpha\in Q_1)}$, and a morphism as
$\left(\begin{smallmatrix}
f_1\\
\vdots\\
f_n
\end{smallmatrix}\right)$. For path $p = \alpha_l\cdots \alpha_1$ with each
$\alpha_i\in Q_1$, let $X_p$ denote the $A$-map $X_{\alpha_l}\cdots
X_{\alpha_1}: X_{s(\alpha_1)}\rightarrow X_{e(\alpha_l)}$.

{\it A
representation $X$ of $(Q, I)$ over $A$} is a representation $X$ of
$Q$ over $A$, such that  $X_{\rho_i}=0$ for each $\rho_i\in \rho$.
Denote by ${\rm rep}(Q, I, A)$ the category of
finite-dimensional representations of $(Q, I)$ over $A$. A morphism
$f$ is a monomorphism (resp., an epimorphism, an isomorphism) if and
only if each $f_i$ is injective (resp., surjective, an isomorphism). A sequence $0\rightarrow X\stackrel
{f}\rightarrow Y\stackrel {g}\rightarrow Z \rightarrow 0$ of
morphisms in ${\rm rep}(Q, I, A)$ is exact if and only if
$0\rightarrow X_i\stackrel {f_i}\rightarrow Y_i\stackrel
{g_i}\rightarrow Z_i \rightarrow 0$ is exact in $A$-mod for each
$i\in Q_0$.

\begin{lem} \label{mod=rep} {\rm([4, p.57], [28, p.44])} \ We have
an equivalence $\m\mbox{-}{\rm mod}\cong {\rm rep}(Q, I, A)$ of
categories. \end{lem}

Throughout we will {\bf identify a $\m$-module with a representation
of $(Q, I)$ over $A$}.

\vskip5pt

As an object of ${\rm rep}(Q, I, k)$ we have $P(v) = (e_i(kQ/I)e_v,
\ P(v)_\alpha,  \ i\in Q_0, \ \alpha\in Q_1),$ where $P(v)_\alpha: \
e_{s(\alpha)}(kQ/I)e_v \rightarrow e_{e(\alpha)}(kQ/I)e_v$ is the
$k$-map sending path $w$ to $\alpha w$. Consider functors \
$-\otimes P(v): A\mbox{-}{\rm mod} \rightarrow {\rm rep}(Q, I, A)$
and $-_v: {\rm rep}(Q, I, A)\rightarrow A\mbox{-}{\rm mod} \ \ \
(\mbox{by taking the $i$-th branch)}.$

\begin{lem} \label{adj} $(1)$ \ $(-\otimes P(v), -_v)$
is an adjoint pair.

\vskip5pt

$(2)$ \   The indecomposable projective $\m$-modules are exactly
$P\otimes P(v)$, where $P$ are indecomposable projective
$A$-modules. In particular, branches of projective $\m$-modules are
projective $A$-modules.
\end{lem}
\noindent{\bf Proof.} \ Let $X=(X_i, \ X_\alpha, \ i\in Q_0, \
\alpha\in Q_1)\in \m${\rm-mod}, $M\in A${\rm-mod}, and $f = (f_i, \
i\in Q_0)\in \Hom_\m(M\otimes P(v), X)$. Then $f_v\in \Hom_A (M,
X_v)$. Since $M\otimes P(v) = (M\otimes e_i(kQ/I)e_v,  \ {\rm
id}_M\otimes P(v)_\alpha, \ i\in Q_0, \ \alpha\in Q_1)$, it follows
from $(1.2)$ that for $i\ne v$ we have
$$f_i = \begin{cases} 0, & \mbox{if there are no non-zero paths from} \ v \ \mbox{to} \ i; \\
m\otimes p\mapsto X_pf_v(m), & \mbox{if there is a non-zero path} \
p \ \mbox{from} \ v \ \mbox{to} \ i.\end{cases} \eqno(1.3)$$  By
$(1.3)$ we see that $f\mapsto f_v$ gives an injective map
$\Hom_\m(M\otimes P(v), X) \rightarrow \Hom_A (M, X_v)$. It is also
surjective, since for a given $f_v\in \Hom_A (M, X_v)$, $f = (f_i, \
i\in Q_0):  M\otimes P(v)\rightarrow X$ is indeed a morphism in
${\rm rep}(Q, I, A)$, where $f_i$ is given by $(1.3)$.

$(2)$ is clear. $\s$

\subsection{} {\it A complete $A$-projective
resolution} is an exact sequence of finitely generated projective
$A$-modules $P^\bullet = \cdots \rightarrow P^{-1}\rightarrow P^{0}
\stackrel{d^0}{\rightarrow} P^{1}\rightarrow \cdots,$ such that
${\rm Hom}_A(P^\bullet, A)$ is also exact. An $A$-module $M$ is {\it
Gorenstein-projective}, if there is a complete $A$-projective
resolution $P^\bullet$ such that $M\cong \operatorname{Ker}d^0$
([13]).  Let $\mathcal {P}(A)$ be the full subcategory of $A$-mod of
projective modules, and $\mathcal {GP}(A)$ the full subcategory of
$A$-mod consisting of Gorenstein-projective modules. Then $\mathcal
{GP}(A)\subseteq \ ^\perp A: = \{X\in A\mbox{-}{\rm mod} \ | \ {\rm
Ext}^i_A(X, A) = 0, \ \forall \ i\ge 1\}$; and $\mathcal {GP}(A) =
A$-mod if and only if $A$ is self-injective. If ${\rm gl.dim} A <
\infty$ then $\mathcal {GP}(A) = \mathcal {P}(A)$ (but the converse
is {\it not} true); and if $A$ is {\it a Gorenstein algebra} (i.e.,
${\rm inj.dim}\ _AA< \infty$ and ${\rm inj.dim} \ A_A < \infty$)
then \ $\mathcal {GP}(A)$ is contravariantly finite in $A$-mod and $\mathcal {GP}(A) = \ ^\perp A$ ([14, 11.5.4, 11.5.3]) (but the
converse is {\it not} true). Note that $\mathcal {GP}(A)$ is a
resolving subcategory of $A$-mod, i.e., $\mathcal
{GP}(A)\supseteq\mathcal {P}(A)$, $\mathcal {GP}(A)$ is closed under
extensions, the kernels of epimorphisms, and direct summands ([17]),
and that $\mathcal {GP}(A)$ is a Frobenius category with relative
projective-injective objects being projective $A$-modules ([6]).

\section {\bf Monic monomial representations}

\begin{defn} \label{maindef} \ A representation $X = (X_i, \
X_{\alpha}, \ i\in Q_0, \ \alpha\in Q_1)$ of $(Q, I)$ over $A$ is a
monic {\rm (}monomial{\rm)} representation, or a monic {\rm
(}monomial{\rm)} $\Lambda$-module, provided that $X$ satisfies the
conditions$:$

\vskip5pt

${\rm (m1)}$  \ For each $i\in Q_0$,  the sum $\sum\limits_{\begin
{smallmatrix} \alpha\in \ha(\to i)
\end{smallmatrix}}\Ima
X_\alpha$ is a direct sum $\bigoplus\limits_{\begin {smallmatrix}
\alpha\in \ha(\to i) \end{smallmatrix}}\Ima X_\alpha;$

\vskip5pt

${\rm (m2)}$  \  For each $\alpha\in Q_1$, there holds $\Ker
X_\alpha =\sum\limits_{q\in K_\alpha} \Ima X_q,$ where $K_\alpha$ is
as in $(2.1)$.
\end{defn}

Easy to see that  ${\rm (m1)}$ and ${\rm (m2)}$ are independent.  If
$K_\alpha = \emptyset$, i.e., $\alpha$ is not the last arrow of all
$\rho_i\in \rho$,  then $\sum\limits_{\begin {smallmatrix} q\in
K_\alpha\end{smallmatrix}} \Ima X_q$ is understood to be zero. We
keep this convention throughout. In this case  ${\rm (m2)}$ says
that $X_\alpha$ is injective $($this is the case if $s(\alpha)$ is a
source$)$. In general, ${\rm (m2)}$ exactly says that \
$\bigoplus\limits_{q\in K_\alpha} X_{s(q)}\stackrel{(X_q)_{q\in
K_\alpha}}\longrightarrow X_{s(\alpha)}\stackrel{X_\alpha}
\longrightarrow X_{e(\alpha)}$ \ is an  exact sequence.

\vskip5pt

Denote by ${\rm mon}(Q, I, A)$ the full subcategory of ${\rm rep}(Q,
I, A)$ consisting of monic representations, which is called {\it the
monomorphism category of $(Q, I)$ over $A$}. Since ${\rm mon}(Q, I,
A)$ is closed under the direct summands, it is a Krull-Schmidt
category ([28]).

\begin{exm} \label{exmmonic} \ $(1)$ \ If $I=0$, then ${\rm(m1)}$ and ${\rm(m2)}$ exactly says
that $(X_{\alpha})_{\alpha\in \ha(\to i)}: \
\bigoplus\limits_{\alpha\in \ha(\to i)} X_{s(\alpha)} \rightarrow
X_i$ is injective for each $i\in Q_0$ {\rm([25])}.

\vskip5pt

$(2)$ \ Indecomposable projective $kQ/I$-module $P(v) =
(e_i(kQ/I)e_v, \ P(v)_\alpha,  \ i\in Q_0, \ \alpha\in Q_1)\in {\rm
rep}(Q, I, k)$ is monic for each $v\in Q_0$, and this gives all the
indecomposable monic $kQ/I$-modules $($see Corollary
{\rm\ref{A=k}}$)$.

\vskip5pt

$(3)$ \ Let $M$ be an arbitrary $A$-module. Then $M\otimes P(v) =
(M\otimes e_i(kQ/I)e_v,  \ {\rm Id}_M\otimes P(v)_\alpha, \ i\in
Q_0, \ \alpha\in Q_1)$ is a monic $\m$-module for each $v\in Q_0$.
Thus  projective $\m$-modules are monic $($cf. Lemma
{\rm\ref{adj}}$(2))$. So, monic $\m$-modules are much more than
projective $\m$-modules, if $A$ is not semisimple.

\vskip5pt

$(4)$  \ Let $A : = k[x]/\langle x^2\rangle$, and $\m: = A\otimes
kQ/I$ with
\[\xymatrix {Q: = \ \ \ 4\ar[r]^-{\gamma}& 3 \ar @/^/[rr]^-{\beta_1} \ar @/_/[rr]_-{\beta_2} &&
2 \ar[r]^-{\alpha}& 1}\] and $I: = \langle \beta_1\gamma,
\alpha\beta_2\gamma\rangle$. Note that $A$ has only two
indecomposable modules $A$ and $k: = k\bar x$ {\rm(}up to
ismomorphisms{\rm)}, with $\Hom_A(A, k) = k\bar x, \ \ \Hom_A(k, A)
= k\sigma, \ \ \Hom_A(A, A) = k{\rm Id}_A\oplus k\bar x,$ \ where $0
\rightarrow k\stackrel \sigma\rightarrow A \stackrel {\bar x}
\rightarrow k \rightarrow 0$ is the canonical exact sequence. Then
\[\xymatrix {X: = \ \ \ k\ar[rr]^-{X_\gamma = \sigma }&& A \ar
@/^/[rr]^-{X_{\beta_1}=\binom{0}{\bar x}} \ar[rr]_-{X_{\beta_2} =
\binom{\rm Id}{\bar x}}
 & & A\oplus k\ar[rr]^-{X_{\alpha}=\left(\begin{smallmatrix}\bar x& 0\\0&{\rm Id}\end{smallmatrix}\right)} && k\oplus k}\]
is a monic $\m$-module.
\end{exm}

The following result give properties of the structural $A$-maps of
monic monomial representations.

\begin{thm}\label{m4} \  Let  $X = (X_i, \ X_{\alpha}, \ i\in Q_0, \ \alpha\in Q_1)$ be
a monic $\Lambda$-module, and $p$ a non-zero path with $l(p) \ge 1$.
Then

\vskip5pt

$(1)$ \ \  $\Ker X_p = \sum\limits_{q \in K_p}\Ima X_q,$ where $K_p$
is as in $(1.1)$. Thus, if $s(p)$ is a source then $X_p$ is
injective.

\vskip5pt

$(2)$ \ \  $\Ker X_p = (\bigoplus\limits_{\beta\in B_1}\Ima X_\beta)
\oplus (\bigoplus\limits_{\beta\in B_2}X_\beta(\Ker X_{p\beta})),$
where\begin{align*}B_1: & = \{\beta \in \mathcal{A}(\to s(p)) \ | \
\beta = {\rm la}(q) \ \mbox{for some} \ q\in \mathcal{P}( \to
s(p)), \ p \beta\in I \}, \\
B_2: & = \{\beta \in \mathcal{A}(\to s(p)) \ | \ \beta = {\rm la}(q)
\ \mbox{for some} \ q\in \mathcal{P}( \to s(p)), \ p\beta\notin I, \
p q\in I\}\end{align*}

{\rm(}If $B_1 = \emptyset = B_2$, then it says that $X_p$ is
injective{\rm)}.

\vskip5pt

$(3)$ \ \ ¡¡Let  $j$ and $i$ be distinct vertices of $Q$. Then
$\sum\limits_{q\in \mathcal{P}(j\to i)}\Ima X_q =
\bigoplus\limits_{q\in \hp( j\to i)}\Ima X_q.$
\end{thm}

\noindent{\bf Proof.} $(1)$ \ Use induction on the length $l(p)$. If
$l(p)=1$, then the assertion is ${\rm(m2)}$. Suppose that $l(p)\ge
2$ with $p = p'\alpha$ and $\alpha\in Q_1$. Clearly $\sum\limits_{q
\in K_p}\Ima X_q \subseteq\Ker X_p$. Let $x\in \Ker X_p$. Then
$X_p(x)=X_{p'}X_\alpha (x)=0,$ i.e., $X_\alpha (x)\in \Ker X_{p'}$.
Since $l(p')<l(p)$, by the induction $X_\alpha (x)\in\Ker X_{p'} =
\sum\limits_{q \in K_{p'}}\Ima X_q.$ For each path $q\in K_{p'}: =
\{q\in\mathcal{P}( \to s(p')) \ | \ p'q\in I\}$, we consider the
last arrow $\beta$ of $q$ with $q = \beta q'$, where $q'$ is a path,
possibly of length $0$. Graphically we have
\[
\xymatrix @R=1.7pc  @C=2pc {
& & s(\alpha) = s(p)\ar[d]^-\alpha &  \\
s(q') \ar @{~>}[r]^-{q'} & s(\beta) \ar[r]^-\beta & e(\alpha) =
s(p')\ar@{~>}[r]^-{p'}&e(p).}\] We divide all these paths $q$ into
two classes via $p'\beta\in I$, or $p'\beta\notin I$. Put
\begin{align*}B_1: & = \{\beta \in \mathcal{A}(\to s(p')) \ | \ \beta = {\rm la}(q) \ \mbox{for some} \ q\in \mathcal{P}( \to
s(p')), \ p' \beta\in I \}, \\
B_2: & = \{\beta \in \mathcal{A}(\to s(p')) \ | \ \beta = {\rm
la}(q) \ \mbox{for some} \ q\in \mathcal{P}( \to s(p')), \
p'\beta\notin I, \ p' q\in I\}.\end{align*} Then $B_1\cap B_2 =
\emptyset$. For those paths $q = \beta q' \in K_{p'}$ such that
$\beta\in B_1$ (this contains the case  $l(q') = 0$), we use the
inclusion $\Ima X_q \subseteq \Ima X_\beta$. Since for $\beta\in
B_1$ we have $p'\beta \in I$, it follows that

$$\sum\limits_{\begin{smallmatrix} q \in K_{p'}\\ {\rm la}(q)\in B_1 \end{smallmatrix} }\Ima
X_q = \sum\limits_{\beta\in B_1}\Ima X_\beta \stackrel{\rm(m1)} =
\bigoplus\limits_{\beta\in B_1}\Ima X_\beta.$$ Then by {\rm(m1)} we
have

$$X_\alpha(x)\in\Ker
X_{p'} = (\bigoplus\limits_{\beta\in B_1}\Ima X_\beta) \oplus
(\bigoplus\limits_{\beta\in B_2}X_\beta
(\sum\limits_{\begin{smallmatrix} q'\in \mathcal{P}( \to s(\beta))\\
p'\beta q'\in I
\end{smallmatrix} } \Ima X_{q'})).\eqno(*)$$

\vskip5pt

Since  $p' \alpha = p \notin I$, we have $\alpha\notin B_1$. It
follows from $(*)$ that either $X_\alpha(x) = 0$ (if $\alpha\notin
B_2$), or (if $\alpha\in B_2$)
$$X_\alpha(x)\in
X_\alpha (\sum\limits_{\begin{smallmatrix} q'\in \mathcal{P}(\to
s(\alpha))\\  p'\alpha q'\in I
\end{smallmatrix}} \Ima X_{q'}).$$
In the both cases there are $x_{q'}\in X_{s(q')}$ such that

$$x- \sum\limits_{\begin{smallmatrix} q' \in \mathcal{P}( \to s(\alpha))\\ p q'\in I \end{smallmatrix} } X_{q'}(x_{q'})\in \Ker X_\alpha.$$
Since  $\Ker X_\alpha =\sum\limits_{\begin{smallmatrix} q \in
\mathcal{P}(\to s(\alpha))\\ \alpha q\in I \end{smallmatrix} }\Ima
X_{q}$ by
${\rm(m2)}$,  it follows that $x\in \sum\limits_{\begin{smallmatrix} q \in \mathcal{P}( \to s(p))\\
pq\in I \end{smallmatrix} }\Ima X_q$, i.e., the assertion $(1)$
holds.

\vskip5pt

$(2)$ in fact has been proved above (see the equality in $(*)$).

\vskip5pt

$(3)$ \ We first prove the following

\vskip5pt

{\bf Claim}: \  If
$\sum\limits_{\begin{smallmatrix} q\in \mathcal P(j\to s(p))\\
pq\notin I\end{smallmatrix}}X_q(x_q) \in \Ker X_p$ with each $x_q\in
X_{s(q)}$, then $X_pX_q(x_q)= 0$ for each $q$.

\vskip5pt

Set $l_p: = 0$ if $s(p)$ is a source; otherwise $l_p$ is defined to
be the maximal length of all the paths $q$ of length $\ge 1$ to
$s(p)$. We stress that  $l_p \ne {\rm max}\{ l(q) \ | \ q\in
\hp(\rightarrow s(p))\}$ in general, since $q\in \hp(\rightarrow
s(p))$ means that $q$ is a non-zero path (considering all the paths
rather than only non-zero paths here will make the argument below
easier). Use induction on $l_p$. If $l_p = 0$ then the assertion
trivially holds.

\vskip5pt

Suppose that $l_p \ge 1$. Put $x: = \sum\limits_{\begin{smallmatrix}
q\in \mathcal P(j\to s(p))\\ pq\notin I\end{smallmatrix}}X_q(x_q)$.
Note that the index set $\{q\in \mathcal P(j\to s(p)) \ | \ pq\notin
I\}$ is a disjoint union of
$$C_1:=\{\beta\in \mathcal A(j\to s(p)) \ | \ p\beta\notin I\}$$ and $$\{q\in \mathcal P(j\to s(p)) \ | \ pq\notin I, \ l(q)\ge 2\}.$$
Denote by $C_2$ the set of the last arrows of $q\in \{q\in \mathcal
P(j\to s(p)) \ | \ pq\notin I, \ l(q)\ge 2\}.$ Since any arrow
$\beta$ in $C_2$ does not starts from $j$, it follows that $C_1\cap
C_2 = \emptyset$. Thus $x = x_1+x_2$ with $$x_1: =
\sum\limits_{\beta\in C_1}X_\beta(x_\beta), \ \ \ \ x_2 =
\sum\limits_{\beta\in C_2}X_\beta(\sum\limits_{\begin{smallmatrix}
q\in \mathcal P(j\to s(\beta))\\ p\beta q\notin
I\end{smallmatrix}}X_q(x_q)).$$

\vskip5pt

On the other hand, by $(2)$ we have $\Ker X_p =
(\bigoplus\limits_{\beta\in B_1}\Ima X_\beta) \oplus
(\bigoplus\limits_{\beta\in B_2}X_\beta(\Ker X_{p\beta})),$
where\begin{align*}B_1: & = \{\beta \in \mathcal{A}(\to s(p)) \ | \
\beta = {\rm la}(q) \ \mbox{for some} \ q\in \mathcal{P}( \to
s(p)), \ p \beta\in I \}, \\
B_2: & = \{\beta \in \mathcal{A}(\to s(p)) \ | \ \beta = {\rm la}(q)
\ \mbox{for some} \ q\in \mathcal{P}( \to s(p)), \ p\beta\notin I, \
p q\in I\}\end{align*} Thus $x = y_1 + y_2$ with $y_1\in
\sum\limits_{\beta\in B_1}\Ima X_\beta$ and $y_2\in
\sum\limits_{\beta\in B_2}X_\beta(\Ker X_{p\beta}).$ Note that
$$C_1\cap C_2 = \emptyset, \ \  B_1\cap B_2 = \emptyset, \ \ B_1\cap C_1 = \emptyset, \ \ B_1\cap C_2 = \emptyset.$$
By $x_1+x_2 = y_1+y_2$ and ${\rm (m1)}$ we see that $y_1 = 0$. So $x
= x_1+x_2 \in \sum\limits_{\beta\in B_2}X_\beta(\Ker X_{p\beta})$.
We discuss each summand $X_q(x_q)$ of $$x =
\sum\limits_{\begin{smallmatrix} q\in \mathcal P(j\to s(p))\\
pq\notin I\end{smallmatrix}}X_q(x_q) =  \sum\limits_{\beta\in
C_1}X_\beta(x_\beta) + \sum\limits_{\beta\in
C_2}X_\beta(\sum\limits_{\begin{smallmatrix} w\in \mathcal P(j\to
s(\beta))\\ p\beta w\notin I\end{smallmatrix}}X_w(x_w))$$ in several
cases. Our aim is to prove that $X_pX_q(x_q) = 0$ for each summand
$X_q(x_q)$ of $x$.

\vskip5pt

Case 1. \ If $\beta\in C_1\setminus B_2$, then by ${\rm (m1)}$ we
see that the summand $X_\beta(x_\beta)$ of $x$ is $0$, and hence
$X_pX_\beta(x_\beta) = 0.$

\vskip5pt

Case 2. \ If $\beta\in C_1\cap B_2$, then by ${\rm (m1)}$ we see
that $X_\beta(x_\beta)= X_\beta(x'_\beta)$ with $x'_\beta\in \Ker
X_{p\beta}$, it follows that
$$X_pX_\beta(x_\beta) = X_pX_\beta(x'_\beta) = X_{p\beta}(x'_\beta) = 0.$$

\vskip5pt

Case 3. \ If $\beta\in C_2\setminus B_2$, then by ${\rm (m1)}$ we
see that the summand $X_\beta(\sum\limits_{\begin{smallmatrix} w\in
\mathcal P(j\to s(\beta))\\ p\beta w\notin
I\end{smallmatrix}}X_w(x_w))$ of $x$ is $0$, and hence
$$\sum\limits_{\begin{smallmatrix} w\in \mathcal P(j\to s(\beta))\\
p\beta w\notin I\end{smallmatrix}}X_w(x_w) \in \Ker X_\beta
\subseteq \Ker X_{p\beta}.$$ Since $l_{p\beta} < l_p$, by induction
 $X_{p\beta}X_w(x_w)= 0$, i.e., $X_pX_\beta X_w(x_w)= 0$.
\vskip5pt

Case 4. \ If $\beta\in C_2\cap B_2$, then by ${\rm (m1)}$ we see
that $X_\beta(\sum\limits_{\begin{smallmatrix} w\in \mathcal P(j\to
s(\beta))\\ p\beta w\notin I\end{smallmatrix}}X_w(x_w))=
X_\beta(x'_\beta)$ with $x'_\beta\in \Ker X_{p\beta}$. It follows
that
 $$X_pX_\beta(\sum\limits_{\begin{smallmatrix} w\in \mathcal P(j\to s(\beta))\\ p\beta w\notin I\end{smallmatrix}}X_w(x_w)) =
 X_pX_\beta(x'_\beta)= X_{p\beta}(x'_\beta) = 0,$$
 i.e., $\sum\limits_{\begin{smallmatrix} w\in \mathcal P(j\to s(\beta))\\ p\beta w\notin I\end{smallmatrix}}X_w(x_w)\in \Ker X_{p\beta}$.
Since $l_{p\beta} < l_p$, by induction $X_{p\beta}X_w(x_w)= 0$,
i.e., $X_pX_\beta X_w(x_w)= 0$.

\vskip5pt

All together {\bf Claim} is proved.

\vskip5pt

Now, let $\sum\limits_{p\in \hp(j\rightarrow i)} X_p(x_{p})=0$ with
each $x_{p}\in X_j$. Then we have
\begin{align*} 0 &= \sum\limits_{p\in \hp(j\rightarrow i)} X_p(x_{p})
=\sum\limits_{\alpha \in  \ha(j\rightarrow i)} X_\alpha(x_{\alpha})
+\sum\limits_{p \in \hp (j\rightarrow i)\setminus \ha(j\to i)}
X_p(x_{p}) \\ & =\sum\limits_{\alpha \in \ha(j\rightarrow i)}
X_\alpha(x_{\alpha})+ \sum\limits_{\beta \in \ha(s(\beta)\to i)}
X_\beta(\sum \limits_{\begin{smallmatrix}q\in \hp(j\rightarrow
s(\beta))\\
\beta q\notin I\end{smallmatrix}}X_q(x_{\beta q})).\end{align*}
Since all those arrows $\beta$ do not start from $j$, by ${\rm(m1)}$
we have $X_\alpha(x_{\alpha})=0$ for $\alpha \in \ha(j\rightarrow
i)$ and $X_\beta(\sum \limits_{\begin{smallmatrix}q\in
\hp(j\rightarrow
s(\beta))\\
\beta q\notin I\end{smallmatrix}}X_q(x_{\beta q}))=0$ for all the
$\beta$'s above. Thus $\sum \limits_{\begin{smallmatrix}q\in
\hp(j\rightarrow
s(\beta))\\
\beta q\notin I\end{smallmatrix}}X_q(x_{\beta q})\in {\rm Ker}
X_\beta.$ It follows from {\bf Claim} that each $X_\beta
X_q(x_{\beta q})=0$. This proves $(3)$. $\s$

\section {\bf Resolvability}

A full subcategory $\mathcal X$ of $\m$-mod is {\it resolving} if
$\mathcal X$ contains all the projective $\m$-modules, $\mathcal X$
is closed under extensions, the kernels of epimorphisms, and direct
summands ([2]).

\begin{thm} \label{resolving} \ The monomorphism category  ${\rm mon}(Q, I, A)$ is a resolving
subcategory of $\m$-mod.

\vskip5pt

Thus, ${\rm mon}(Q, I, A)$ has enough projective objects, which are
exactly projective $\m$-modules. In particular, each branch of
projective objects of ${\rm mon}(Q, I, A)$ is a projective
$A$-module.
\end{thm} \noindent{\bf Proof.}  By
Example \ref{exmmonic}$(3)$ ${\rm mon}(Q, I, A)$ contains $\mathcal
P(\m)$. Clearly ${\rm mon}(Q, I, A)$ is closed under direct
summands.  Let $0\rightarrow X \stackrel {f }\rightarrow Y \stackrel
{g}\rightarrow Z \rightarrow 0$ be an exact sequence in ${\rm
rep}(Q, I, A)$ with $X, Z\in{\rm mon}(Q, I, A)$.

\vskip5pt

{\bf Claim 1}:  \ $Y$ satisfies ${\rm (m1)}$, i.e.,
$\sum\limits_{\begin {smallmatrix} \alpha\in \ha(\to i)
\end{smallmatrix}}\Ima
Y_\alpha = \bigoplus\limits_{\begin {smallmatrix} \alpha\in \ha(\to
i) \end{smallmatrix}}\Ima Y_\alpha$ for each $i\in Q_0$.

\vskip5pt

The argument is a diagram chasing. By $(1.2)$  we have a commutative
diagram with exact rows
\[\xymatrix {0 \ar[r] &  \bigoplus\limits_{\begin {smallmatrix} \alpha\in
\ha(\to i)
\end{smallmatrix}}
X_{s(\alpha)} \ar[d]_-{(X_{\alpha})_{\alpha \in \ha(\to
i)}}\ar[r]^-{\oplus f_{s(\alpha)}}& \bigoplus\limits_{\begin
{smallmatrix} \alpha\in \ha(\to i)
\end{smallmatrix}}
Y_{s(\alpha)}\ar[d]^-{(Y_{\alpha})_{\alpha \in \ha(\to
i)}}\ar[r]^-{\oplus g_{s(\alpha)}}&\bigoplus\limits_{\begin
{smallmatrix} \alpha\in \ha(\to i)
\end{smallmatrix}}
Z_{s(\alpha)} \ar[d]^-{(Z_{\alpha})_{\alpha \in \ha(\to i)}}\ar[r] &
0 \\   0 \ar[r] & X_i\ar[r]^-{f_i} & Y_i\ar[r]^{g_i}& Z_i\ar[r]&
0.}\eqno(3.1)\]

\noindent Note that $\sum\limits_{\begin {smallmatrix} \alpha\in
\ha(\to i)
\end{smallmatrix}}\Ima
Y_\alpha = \Ima (Y_\alpha)_{\alpha\in \ha(\to i)}$. Let
$\sum\limits_{\begin {smallmatrix} \alpha\in \ha(\to i)
\end{smallmatrix}}Y_\alpha(y_\alpha)=0 $ for some $y_\alpha\in Y_{s(\alpha)}.$
By the right square above we have
$$\sum\limits_{\begin {smallmatrix} \alpha\in
\ha(\to i)
\end{smallmatrix}}g_iY_\alpha(y_\alpha)=\sum\limits_{\begin {smallmatrix} \alpha\in
\ha(\to i)
\end{smallmatrix}} Z_\alpha g_{s(\alpha)}(y_\alpha) =0.$$
 Since $Z$ satisfies ${\rm(m1)}$, we have $Z_\alpha g_{s(\alpha)}(y_\alpha)=0$ for each $\alpha$. By ${\rm(m2)}$ on $Z$ we
 get
$g_{s(\alpha)}(y_\alpha)=\sum\limits_{q\in K_\alpha} Z_q(z_q)$ for
some $z_q\in Z_{s(q)}$. Since $g_{s(q)}: Y_{s(q)}\rightarrow
Z_{s(q)}$ is surjective, there exists $y_q\in Y_{s(q)}$ such that
$z_q=g_{s(q)}(y_q)$. So by $(1.2)$ we have (note that $e(q) =
s(\alpha)$)
$$g_{s(\alpha)}(y_\alpha)=\sum\limits_{q\in K_\alpha} Z_q(g_{s(q)}(y_q)) =\sum\limits_{q\in K_\alpha}  g_{s(\alpha)}(Y_q(y_q))$$
for each $\alpha\in \ha(\to i)$, i.e., \ $y_\alpha-\sum\limits_{q\in
K_\alpha} Y_q(y_q)\in \Ker g_{s(\alpha)} = \Ima f_{s(\alpha)}.$
Hence
$$y_\alpha-\sum\limits_{q\in K_\alpha} Y_q(y_q) =
 f_{s(\alpha)}(x_\alpha)$$ for some $x_\alpha\in X_{s(\alpha)}$.
Applying $Y_\alpha$ and by $Y_\alpha Y_q=0$ for each $q\in
K_\alpha$, we have
$$Y_{\alpha}(y_\alpha)= Y_{\alpha}(
f_{s(\alpha)}(x_\alpha))=f_i(X_\alpha(x_\alpha)).$$ By
$\sum\limits_{\begin {smallmatrix} \alpha\in \ha(\to i)
\end{smallmatrix}}Y_\alpha(y_\alpha)=0$  we have $\sum\limits_{\begin {smallmatrix} \alpha\in
\ha(\to i)
\end{smallmatrix}}f_i(X_\alpha(x_\alpha)) =0,$ and hence $\sum\limits_{\begin {smallmatrix} \alpha\in
\ha(\to i)
\end{smallmatrix}}X_\alpha(x_\alpha)=0$ since $f_i$ is injective.  Since $X$ satisfies ${\rm(m1)}$, we have  $X_\alpha(x_\alpha)=0$
for each $\alpha\in \ha(\to i)$, and hence
$Y_{\alpha}(y_\alpha)=f_i(X_\alpha(x_\alpha))=0.$ This proves {\bf
Claim 1}.

\vskip5pt

{\bf Claim 2}: \ $Y$ satisfies ${\rm (m2)}$, i.e., $\Ker Y_\alpha
=\sum\limits_{q\in K_\alpha} \Ima Y_{q}$ for each $\alpha\in Q_1$.

\vskip5pt

In fact, by $(1.2)$ we have a commutative diagram with exact rows
\[\xymatrix {& 0\ar[d] & 0\ar[d] & 0\ar[d] & \\
0 \ar[r] &  \bigoplus\limits_{q\in K_\alpha} X_{s(q)}
\ar[d]_-{(X_q)_{q\in K_\alpha}}\ar[r]^-{\oplus f_{s(q)}}&
\bigoplus\limits_{q\in K_\alpha} Y_{s(q)} \ar[d]^-{(Y_q)_{q\in
K_\alpha}}\ar[r]^-{\oplus g_{s(q)}}&\bigoplus\limits_{q\in K_\alpha}
Z_{s(q)} \ar[d]^-{(Z_q)_{q\in K_\alpha}}\ar[r] & 0
\\   0 \ar[r] &
X_{s(\alpha)}\ar[d]_-{X_{\alpha}}\ar[r]^-{f_{s(\alpha)}} &
Y_{s(\alpha)}\ar[d]^-{Y_{\alpha}}\ar[r]^{g_{s(\alpha)}}&
Z_{s(\alpha)}\ar[d]^-{Z_{\alpha}}\ar[r]& 0
\\   0 \ar[r] & X_{e(\alpha)} \ar[r]^-{f_{e(\alpha)}}\ar[d] & Y_{e(\alpha)} \ar[r]^-{g_{e(\alpha)}}\ar[d]& Z_{e(\alpha)}\ar[r]\ar[d] &
0 \\  & 0 & 0 & 0 & }\]

\vskip5pt

\noindent Viewing columns as complexes such that $X_{s(\alpha)}$ is
at position zero, the diagram above gives a short exact sequence of
complexes, and hence it induces a long exact sequence of
cohomologies. The assumption that $X$ and $Z$ satisfy ${\rm(m2)}$
means that the $0$-th cohomologies of the complexes on the left and
on the right are zero. It follows that the $0$-th cohomologies of
the complex at the middle is also zero, i.e., $Y$ satisfies
${\rm(m2)}$. This proves {\bf Claim 2}, and hence that ${\rm mon}(Q,
I, A)$ is closed under extensions.

\vskip5pt

Now we prove that ${\rm mon}(Q, I, A)$ is closed under the kernels
of epimorphisms. Let $0\rightarrow X \stackrel {f }\rightarrow Y
\stackrel {g}\rightarrow Z \rightarrow 0$ be an exact sequence in
${\rm rep}(Q, I, A)$ with $Y, Z\in{\rm mon}(Q, I, A)$. We need to
prove that $X$ satisfies ${\rm (m1)}$ and ${\rm(m2)}$.  Let
$\sum\limits_{\begin {smallmatrix} \alpha\in \ha(\to i)
\end{smallmatrix}}X_\alpha(x_\alpha)=0 $ for some $x_\alpha\in X_{s(\alpha)}.$
By the left square in $(3.1)$ we have $\sum\limits_{\begin
{smallmatrix} \alpha\in \ha(\to i)
\end{smallmatrix}}Y_\alpha f_{s(\alpha)}(x_\alpha)=0.$
Since $Y$ satisfies ${\rm(m1)}$, we have  $Y_\alpha
f_{s(\alpha)}(x_\alpha)=0$ for each $\alpha$. By $(1.2)$ this is
exactly $f_i X_\alpha(x_\alpha)=0$. Since $f_i$ is injective, we
have $X_\alpha(x_\alpha)=0$ for each $\alpha\in \ha(\to i)$. This
proves that $X$ satisfies ${\rm (m1)}$. It remains to prove that $X$
satisfies ${\rm (m2)}$, i.e.,  $\Ker X_\alpha =\sum\limits_{q\in
K_\alpha} \Ima X_{q}$ for each $\alpha\in E$. This follows from the
following sub-lemma.

\vskip5pt

Finally,  by Example \ref{exmmonic}$(3)$  projective $\m$-modules
are monic. From this and the fact that ${\rm mon}(Q, I, A)$ is
closed under the kernels of epimorphisms we see that ${\rm mon}(Q,
I, A)$ has enough projective objects, which are exactly projective
$\m$-modules. $\s$

\begin{sublem} \   Let $0\rightarrow X \stackrel {f }\rightarrow Y \stackrel
{g}\rightarrow Z \rightarrow 0$ be an exact sequence in ${\rm
rep}(Q, I, A)$ with $Y$ and \ $Z$ in ${\rm mon}(Q, I, A)$. Then
$\Ker X_p =\sum\limits_{q\in K_p} \Ima X_{q}$ for each non-zero path
$p$, where $K_p$ is as in $(1.1)$.

\end{sublem}
\noindent{\bf Proof.} Define $l_p$ as in the proof of Theorem
\ref{m4}$(3)$. Use induction on $l_p$. If $l_p: = 0$, then $s(p)$ is
a source and hence $\Ker Y_p = 0$ by Theorem \ref{m4}$(1)$. So $\Ker
X_p = 0$, i.e., the assertion holds.

Suppose that the assertion holds for all the non-zero paths with
$l_p < m \ (m\ge 1)$. We prove the assertion for non-zero path $p$
with $l_p=m.$ Let $x\in \Ker X_p$. By the commutative diagram with
exact rows
\[\xymatrix { 0 \ar[r] & X_{s(p)}\ar[d]^-{X_{p}}\ar[r]^-{f_{s(p)}} & Y_{s(p)}\ar[d]^-{Y_{p}}\ar[r]^{g_{s(p)}}& Z_{s(p)}\ar[d]^-{Z_{p}}\ar[r]&
0 \\   0 \ar[r] & X_{e(p)}\ar[r]^-{f_{e(p)}} &
Y_{e(p)}\ar[r]^{g_{e(p)}}& Z_{e(p)}\ar[r]& 0}\]

\noindent we have $f_{s(p)}(x)\in \Ker Y_p$. By Theorem
\ref{m4}$(2)$ we know that $f_{s(p)}(x)$ is of the form
$$f_{s(p)}(x)=\sum\limits_{\begin {smallmatrix} \beta\in B_1
\end{smallmatrix}}Y_{\beta}(y_\beta)+\sum\limits_{\begin
{smallmatrix} \beta\in B_2
\end{smallmatrix}}Y_{\beta}(y_\beta)$$ where $y_\beta\in
Y_{s(\beta)}$ for each $\beta\in B_1\dot{\cup} B_2$, moreover
$y_\beta\in\Ker Y_{p\beta}$ for each $\beta\in B_2$.  By  $(1.2)$ we
have $$\sum\limits_{\begin {smallmatrix} \beta\in B_1\dot{\cup} B_2
\end{smallmatrix}}Z_{\beta}g_{s(\beta)}(y_\beta)
=\sum\limits_{\begin {smallmatrix} \beta\in B_1\dot{\cup} B_2
\end{smallmatrix}}g_{s(p)}Y_{\beta}(y_\beta) =g_{s(p)}f_{s(p)}(x)=0.$$ Since
$B_1\cap B_2 = \emptyset$ and $Z$ satisfies $(m1)$, we have
$Z_{\beta}g_{s(\beta)}(y_\beta)=0$ for each $\beta\in B_1\dot{\cup}
B_2$. Then by ${\rm(m2)}$ on $Z$ we have
$g_{s(\beta)}(y_\beta)\in\Ker Z_\beta =\sum\limits_{q\in K_\beta}
\Ima Z_q$, so there are some $z_q\in Z_{s(q)}$ such that
$g_{s(\beta)}(y_\beta)=\sum\limits_{q\in K_\beta}Z_q(z_q).$ Since
$g_{s(q)}$ is surjective, there are $y_q\in Y_{s(q)}$ such that $z_q
= g_{s(q)}(y_q).$ Thus
 $$g_{s(\beta)}(y_\beta)=\sum\limits_{q\in K_\beta}
Z_q(g_{s(q)}(y_q))=\sum\limits_{q\in
K_\beta}g_{s(\beta)}(Y_q(y_q)),$$ and hence
$y_\beta-\sum\limits_{q\in K_\beta} Y_q(y_q)\in \Ker g_{s(\beta)}
=\Ima f_{s(\beta)},$ i.e., $y_\beta-\sum\limits_{q\in
K_\beta}Y_q(y_q) = f_{s(\beta)}(x_\beta)$ for some $x_\beta\in
X_{s(\beta)}$. Thus for each $\beta\in B_1\dot{\cup} B_2$ we have $
Y_\beta(y_\beta)=Y_\beta f_{s(\beta)}(x_\beta).$

\vskip5pt

Now for each $\beta\in B_2$, since $y_\beta\in\Ker Y_{p\beta}$, we
have
$$f_{e(p)}X_pX_\beta(x_\beta) = Y_pf_{s(p)} X_\beta(x_\beta) = Y_pY_\beta f_{s(\beta)}(x_\beta) = Y_p Y_\beta(y_\beta) =0,$$
and hence $x_\beta\in \Ker X_{p\beta}$ for each $\beta\in B_2$.
Since $l_{p\beta}<l_p$, by induction $\Ker
X_{p\beta}=\sum\limits_{q\in K_{p\beta}}\Ima X_q$. So
\begin{align*}f_{s(p)}(x)&=\sum\limits_{\begin {smallmatrix}
\beta\in B_1
\end{smallmatrix}} Y_\beta f_{s(\beta)}(x_\beta) +
\sum\limits_{\begin {smallmatrix} \beta\in B_2 \end{smallmatrix}}
Y_\beta f_{s(\beta)}(x_\beta)\\ &=\sum\limits_{\begin {smallmatrix}
\beta\in B_1 \end{smallmatrix}} f_{s(p)}X_\beta(x_\beta) +
\sum\limits_{\begin {smallmatrix} \beta\in B_2
\end{smallmatrix}}f_{s(p)} X_\beta (x_\beta).\end{align*} Since
$f_{s(p)}$ is injective, it follows that
\begin{align*}x & =\sum\limits_{\begin {smallmatrix} \beta\in B_1 \end{smallmatrix}} X_\beta(x_\beta) + \sum\limits_{\begin {smallmatrix} \beta\in B_2 \end{smallmatrix}} X_\beta
(x_\beta)\in  \sum\limits_{\begin {smallmatrix} \beta\in B_1
\end{smallmatrix}} \Ima X_\beta + \sum\limits_{\begin {smallmatrix}
\beta\in B_2 \end{smallmatrix}} X_\beta (\sum\limits_{q\in
K_{p\beta}}\Ima X_q)\\ & \subseteq\sum\limits_{\beta\in B_1}\Ima
X_\beta + \sum\limits_{q\in K_p}\Ima X_q\subseteq \sum\limits_{q\in
K_p}\Ima X_q.\end{align*}

\noindent This completes the proof. $\s$

\vskip5pt  It seems to be more difficult to get all the injective
objects of ${\rm mon}(Q, I, A)$. We can only give some
indecomposable injective objects of ${\rm mon}(Q, I, A)$. This is
needed in the next section.

\begin{prop} \label{injectives} \ Let $N$ be an indecomposable injective $A$-module and  $P(v)$
the indecomposable projective $kQ/I$-module at $v\in Q_0$. Then
$N\otimes P(v)$ is an indecomposable injective object of ${\rm
mon}(Q, I, A)$.
\end{prop}

\noindent{\bf Proof.} \  By Lemma \ref{adj} $N\otimes P(v)$ is an
indecomposable object of ${\rm mon}(Q, I, A)$. Put $L: =
D(A_A)\otimes kQ/I$, where $D: = {\rm Hom}_k(-, k)$. It suffices to
prove that $L$ is an injective object of ${\rm mon}(Q, I, A)$.  Use
induction on $|Q_0|$. Write $L$ as $L = (L_i, \ L_\alpha, \ i\in
Q_0, \ \alpha\in Q_1).$  Let $Q''$ be the quiver by deleting sink
vertex $1$ from $Q$, and $I'': = \langle \rho_i \ | \ e(\rho_i)\neq
1\rangle$. For a $\m$-module $X$, let $X''$ be the $(A\otimes
kQ''/I'')$-module by deleting the first branch $X_1$ from $X$. For a
$\Lambda$-map $f$, we similarly define $f''$.

\vskip5pt

Let $0 \rightarrow X \stackrel{f} \rightarrow Y
\stackrel{g}\rightarrow Z \rightarrow 0$  be an exact sequence and
$h: X\rightarrow L$ a morphism, both in ${\rm mon}(Q, I, A)$. Then
$0 \rightarrow X'' \stackrel{f''} \rightarrow Y''
\stackrel{g''}\rightarrow Z'' \rightarrow 0$ is an exact sequence in
${\rm mon}(Q'', I'', A).$ By induction $L'': = D(A_A)\otimes kQ/I''$
is an injective object of ${\rm mon}(Q'', I'', A)$. So there is a
morphism $u'' = \left(\begin{smallmatrix}u_2\\\vdots\\u_n
\end{smallmatrix}\right): Y''\rightarrow L''$ in ${\rm rep}(Q'',
I'', A)$ such that $h'' = u''f''.$ Thus, it suffices to prove the
following

\vskip5pt

{\bf Claim:} \ there is an $A$-map $u_1: Y_1\rightarrow L_1$ such
that $u = \left(\begin{smallmatrix}u_1\\u_2\\\vdots\\u_n
\end{smallmatrix}\right): Y\rightarrow L$ is a morphism in ${\rm rep}(Q, I, A)$, and that $h_1 = u_1f_1$.

\vskip5pt

Since $L_1$ is a direct of copies of $D(A_A)$, it is an injective
$A$-module. So there is an $A$-map $u_1': Y_1\rightarrow L_1$ such
that $h_1=u_1'f_1$. Consider the $A$-map $(L_\alpha u_{s(\alpha)} -
u_1'Y_\alpha)_{\alpha\in \ha(\to 1)}: \bigoplus\limits_{\alpha\in
\ha(\to 1)} Y_{s(\alpha)} \rightarrow L_1$ and the exact sequence of
$A$-modules $$0 \rightarrow \bigoplus\limits_{\alpha\in \ha(\to 1)}
X_{s(\alpha)} \stackrel{\oplus f_{s(\alpha)}} \longrightarrow
\bigoplus\limits_{\alpha\in \ha(\to 1)} Y_{s(\alpha)}
\stackrel{\oplus g_{s(\alpha)}}\longrightarrow
\bigoplus\limits_{\alpha\in \ha(\to 1)} Z_{s(\alpha)} \rightarrow
0.$$ Since $f$ and $h$ are morphisms, we have
\begin{align*} & (L_\alpha u_{s(\alpha)} - u_1'Y_\alpha)_{\alpha\in \ha(\to 1)} \circ \oplus f_{s(\alpha)}  = (L_\alpha
u_{s(\alpha)} f_{s(\alpha)} - u_1'Y_\alpha f_{s(\alpha)})_{\alpha\in
\ha(\to 1)}\\ & = (L_\alpha u_{s(\alpha)} f_{s(\alpha)} -
u_1'f_1X_\alpha)_{\alpha\in \ha(\to 1)}
 = (L_\alpha h_{s(\alpha)}-
h_1X_\alpha)_{\alpha\in \ha(\to 1)} = 0.\end{align*} It follows that
$(L_\alpha u_{s(\alpha)} - u_1'Y_\alpha)_{\alpha\in \ha(\to 1)}$
factors through $\oplus g_{s(\alpha)}$, i.e., there is  an $A$-map
$v: \bigoplus\limits_{\alpha\in \ha(\to 1)} Z_{s(\alpha)}
\rightarrow L_1$ such that
$$(L_\alpha u_{s(\alpha)} - u_1'Y_\alpha)_{\alpha\in \ha(\to 1)} = v\circ \oplus g_{s(\alpha)}.\eqno(3.2)$$

Put $E$ to be the set of the last arrows of $\rho_i\in \rho$. Write
$\ha(\to 1)$ as a disjoint union $\ha(\to 1) = A_1\dot{\cup} A_2$,
where $A_1: = \ha(\to 1)\setminus E$ and  $A_2: = \ha(\to 1)\cap E$.
For each $\alpha\in A_2$, put $Z^\alpha: = (Z_q)_{q\in K_\alpha}:
\bigoplus\limits_{q\in K_\alpha} Z_{s(q)} \rightarrow
Z_{s(\alpha)}$. Similarly we have $Y^\alpha$ and $L^\alpha$. Clearly
we have
$$Z_\alpha Z^\alpha = 0, \ \ Y_\alpha Y^\alpha = 0, \ \ L_\alpha L^\alpha = 0.\eqno(3.3)$$ By ${\rm (m2)}$
we have $\Ima Z^\alpha = \Ker Z_\alpha, \ \ \ \Ima Y^\alpha = \Ker
Y_\alpha.$

\vskip5pt

Consider the $A$-map $\bigoplus\limits_{\alpha\in A_2} Z^{\alpha}:
\bigoplus\limits_{\alpha\in A_2}(\bigoplus\limits_{q\in K_\alpha}
Z_{s(q)}) \rightarrow \bigoplus\limits_{\alpha\in
A_2}Z_{s(\alpha)}$, also the $A$-map
$$\left(\begin{smallmatrix} 0
\\ \\ \bigoplus\limits_{\alpha\in A_2} Z^\alpha \end{smallmatrix}\right): \bigoplus\limits_{\alpha\in A_2}(\bigoplus\limits_{q\in K_\alpha}
Z_{s(q)})\rightarrow (\bigoplus\limits_{\alpha\in
A_1}Z_{s(\alpha)})\oplus (\bigoplus\limits_{\alpha\in
A_2}Z_{s(\alpha)}) = \bigoplus\limits_{\alpha\in \ha(\to
1)}Z_{s(\alpha)}.$$ Note that
$$\Ima (\left(\begin{smallmatrix} 0
\\ \\ \bigoplus\limits_{\alpha\in A_2} Z^\alpha \end{smallmatrix}\right)) = \bigoplus\limits_{\alpha\in A_2} \Ker
Z_\alpha, \ \ \ \ \ \Ima (\left(\begin{smallmatrix} 0
\\ \\ \bigoplus\limits_{\alpha\in A_2} Y^\alpha \end{smallmatrix}\right)) = \bigoplus\limits_{\alpha\in A_2} \Ker Y_\alpha.$$
Thus $$\left(\begin{smallmatrix}0\\ \\
\bigoplus\limits_{\alpha\in A_2} Z^\alpha \end{smallmatrix}\right) =
\sigma_2\circ \bigoplus\limits_{\alpha\in
A_2}\widetilde{Z^\alpha},\eqno(3.4)$$ where
$\bigoplus\limits_{\alpha\in A_2}\widetilde{Z^\alpha}:
\bigoplus\limits_{\alpha\in A_2}(\bigoplus\limits_{q\in K_\alpha}
Z_{s(q)})\twoheadrightarrow\bigoplus\limits_{\alpha\in A_2} \Ker
Z_\alpha$, and $\sigma_2: \bigoplus\limits_{\alpha\in A_2} \Ker
Z_\alpha\hookrightarrow \bigoplus\limits_{\alpha\in \ha(\to 1)}
Z_{s(\alpha)}$ is the embedding. By $(1.2)$ we have the commutative
diagram
\[\xymatrix {\bigoplus\limits_{\alpha\in A_2}(\bigoplus\limits_{q\in K_\alpha}
Y_{s(q)})\ar[rr]^-{\left(\begin{smallmatrix} 0
\\ \\ \bigoplus\limits_{\alpha\in A_2} Y^\alpha \end{smallmatrix}\right)}\ar[dd]_-{\bigoplus\limits_{\alpha\in
A_2}\bigoplus\limits_{q\in K_\alpha}g_{s(q)}} &&
\bigoplus\limits_{\alpha\in \ha(\to
1)}Y_{s(\alpha)}\ar[dd]^-{\bigoplus\limits_{\alpha\in \ha(\to
1)}g_{s(\alpha)}} \\ \\ \bigoplus\limits_{\alpha\in
A_2}(\bigoplus\limits_{q\in K_\alpha}
Z_{s(q)})\ar[rr]^-{\left(\begin{smallmatrix} 0
\\ \\ \bigoplus\limits_{\alpha\in A_2} Z^\alpha \end{smallmatrix}\right)} &&
\bigoplus\limits_{\alpha\in \ha(\to 1)}Z_{s(\alpha)}}\eqno(3.5)\]
where the vertical maps are surjective. For each $\alpha\in A_2$ we
have the exact sequence
$$0\longrightarrow \Ker Z_\alpha \hookrightarrow
Z_{s(\alpha)}\stackrel{\widetilde{Z_\alpha}}\longrightarrow \Ima
Z_\alpha\longrightarrow 0$$ induced by $Z_\alpha: Z_{s(\alpha)}
\rightarrow Z_1$, and hence we have exact sequence
$$0 \longrightarrow \bigoplus\limits_{\alpha\in A_2} \Ker Z_\alpha \longrightarrow
\bigoplus\limits_{\alpha\in
A_2}Z_{s(\alpha)}\stackrel{\bigoplus\limits_{\alpha\in
A_2}\widetilde{Z_\alpha}} \longrightarrow
\bigoplus\limits_{\alpha\in A_2}\Ima Z_{\alpha}\longrightarrow 0.$$
Thus we have exact sequence
$$0 \longrightarrow \bigoplus\limits_{\alpha\in A_2} \Ker Z_\alpha \stackrel{\sigma_2} \longrightarrow
\bigoplus\limits_{\alpha\in \ha(\to 1)}Z_{s(\alpha)}
\stackrel{\pi}\longrightarrow (\bigoplus\limits_{\alpha\in A_1}
Z_{s(\alpha)})\oplus (\bigoplus\limits_{\alpha\in A_2}\Ima
Z_{\alpha})\longrightarrow 0,$$ where $\pi:
=\left(\begin{smallmatrix} {\rm id}_{\bigoplus\limits_{\alpha\in
A_1} Z_{s(\alpha)}} & 0
\\ 0 & \bigoplus\limits_{\alpha\in
A_2}\widetilde{Z_\alpha} \end{smallmatrix}\right)$.  Composing $v:
\bigoplus\limits_{\alpha\in \ha(\to 1)} Z_{s(\alpha)}
\longrightarrow L_1$ on the right with $\sigma_2\circ
\bigoplus\limits_{\alpha\in A_2} \widetilde{Z^\alpha}\circ
(\bigoplus\limits_{\alpha\in A_2}\bigoplus\limits_{q\in
K_\alpha}g_{s(q)}),$ by $(3.4)$, $(3.5)$, $(3.2)$ and $(3.3)$ we
have
\begin{align*} & v\circ \sigma_2\circ \bigoplus\limits_{\alpha\in A_2} \widetilde{Z^\alpha}\circ (\bigoplus\limits_{\alpha\in
A_2}\bigoplus\limits_{q\in K_\alpha}g_{s(q)})
 =v\circ\left(\begin{smallmatrix}0\\ \\ \bigoplus\limits_{\alpha\in A_2}
Z^\alpha \end{smallmatrix}\right) \circ (\bigoplus\limits_{\alpha\in
A_2}\bigoplus\limits_{q\in K_\alpha}g_{s(q)})
\\ & = v\circ (\bigoplus\limits_{\alpha\in
\ha(\to 1)}g_{s(\alpha)})\circ \left(\begin{smallmatrix}0\\ \\
\bigoplus\limits_{\alpha\in A_2} Y^\alpha \end{smallmatrix}\right)
 = (L_\alpha u_{s(\alpha)}-u_1' Y_{\alpha})_{\alpha\in \ha(\to
1)}\circ \left(\begin{smallmatrix}0\\ \\
\bigoplus\limits_{\alpha\in A_2} Y^\alpha \end{smallmatrix}\right)
\\&=(L_\alpha u_{s(\alpha)}-u_1' Y_{\alpha})_{\alpha\in
A_2}\circ\bigoplus\limits_{\alpha\in A_2} Y^\alpha = (L_\alpha
u_{s(\alpha)}Y^{\alpha}-u_1' Y_\alpha Y^\alpha)_{\alpha\in A_2}\\& =
(L_\alpha u_{s(\alpha)}Y^{\alpha})_{\alpha\in A_2}
 = (L_{\alpha}L^{\alpha}\circ \bigoplus\limits_{q\in K_\alpha}u_{s(q)})_{\alpha\in A_2}
  = 0,\end{align*}
where we have used the following commutative diagram
\[\xymatrix {\bigoplus\limits_{q\in K_\alpha}
Y_{s(q)}\ar[r]^-{Y^\alpha}\ar[d]_-
{\bigoplus\limits_{q\in K_\alpha}u_{s(q)}} & Y_{s(\alpha)}\ar[d]^-{u_{s(\alpha)}} \\
\bigoplus\limits_{q\in K_\alpha} L_{s(q)}\ar[r]^-{L^\alpha} &
L_{s(\alpha)}.}\]  Since $\bigoplus\limits_{\alpha\in A_2}
\widetilde{Z^\alpha}\circ\bigoplus\limits_{\alpha\in
A_2}\bigoplus\limits_{q\in K_\alpha}g_{s(q)}$ is surjective, we get
$v\circ \sigma_2 =0$, thus $v$ factors through $\pi$. Hence there is
an $A$-map $$t: (\bigoplus\limits_{\alpha\in A_1}
Z_{s(\alpha)})\oplus (\bigoplus\limits_{\alpha\in A_2}\Ima
Z_{\alpha})\longrightarrow L_1$$ such that $v = t\pi$. Since $L_1$
is an injective $A$-module, there is an $A$-map $ w:
Z_1\longrightarrow L_1$ such that the diagram
\[\xymatrix{(\bigoplus\limits_{\alpha\in A_1}
Z_{s(\alpha)})\oplus (\bigoplus\limits_{\alpha\in A_2}\Ima
Z_{\alpha}) \ar@{^{(}->}[rr]^-{((Z_\alpha)_{\alpha\in A_1},
\sigma)}\ar [d]_-{t} &&
Z_1 \ar @{-->} [lld]^-{w}\\
L_1  & &}\] commutes, where $\sigma$ is the embedding (note that
$(Z_\alpha)_{\alpha\in A_1}$ is injective by ${\rm (m2)}$, and hence
$((Z_\alpha)_{\alpha\in A_1}, \sigma)$ injective by ${\rm (m1)}$).
All together we have
\begin{align*}\ \ \ \ \ \ \ \ \ \ \ \ \ \ \ \ \ \ \ \ \ \ \ \ \ v & =
t\pi=w((Z_\alpha)_{\alpha\in A_1}, \sigma) \pi =
w((Z_\alpha)_{\alpha\in A_1}, \sigma)
\left(\begin{smallmatrix}{\rm id} & 0 \\
0 & \bigoplus\limits_{\alpha\in A_2}\widetilde{ Z_\alpha}
\end{smallmatrix}\right) \\ & = w ((Z_\alpha)_{\alpha\in A_1},  (Z_\alpha)_{\alpha\in A_2}) = w (Z_\alpha)_{\alpha\in \ha(\to 1)}.  \ \ \ \ \ \ \ \ \ \ \ \ \ \ \ \ \ \ \ \ \ \ \ \ \ \ \ \ \ \ \ \ \ \ \ \ \ \ \ \ \ \ (3.6)\end{align*}

\vskip5pt

Now put $u_1:=u_1'+wg_1: Y_1 \rightarrow L_1$, where $g_1:
Y_1\rightarrow Z_1$ is the $1$-st branch of $g: Y \rightarrow Z$.
Then $u_1f_1=u_1'f_1+wg_1f_1=u_1'f_1=h_1$. It
remains to prove that $u=\left(\begin{smallmatrix} u_1 \\ u_2 \\
\vdots \\ u_n
\end{smallmatrix} \right): Y\rightarrow L $ is a $\Lambda$-map,
i.e., for each arrow $\alpha: s(\alpha) \rightarrow 1$ the diagram
\[\xymatrix{Y_{s(\alpha)}\ar[r]^-{u_{s(\alpha)}}\ar[d]_-
{Y_\alpha} & L_{s(\alpha)} \ar[d]^-{L_\alpha} \\
L_1\ar[r]^-{u_1} & L_1}\] commutes. In fact by $(3.2)$ and $(3.6)$
we have
\begin{align*} &(L_{\alpha}u_{s(\alpha)}-u_1' Y_{\alpha})_{\alpha\in \ha(\to 1)}= v\circ \oplus g_{s(\alpha)} \\ &= w \circ (Z_{\alpha})_{\alpha\in \ha(\to 1)}\circ \oplus g_{s(\alpha)}
 =w \circ ( Z_{\alpha}g_{s(\alpha)})_{\alpha\in \ha(\to 1)}
\\ &= w \circ (g_1Y_{\alpha})_{\alpha\in \ha(\to 1)}
 = (wg_1Y_{\alpha})_{\alpha\in \ha(\to 1)}.\end{align*}
This means that $L_{\alpha}u_{s(\alpha)}-u_1' Y_{\alpha}
=wg_1Y_{\alpha}$, i.e., $L_{\alpha}u_{s(\alpha)} = u_1Y_{\alpha}.$

\vskip5pt

This proves {\bf Claim}, and hence completes the proof. $\s$

\section {\bf Main result}

Let $A$ be a finite-dimensional algebra over field $k$, $Q$ a finite
acyclic quiver, $I$ an ideal of $kQ$ generated by monomial
relations, and $\m: =A\otimes kQ/I$. Note that $\m$ is not
necessarily a Gorenstein algebra. Our aim is to characterize the
Gorenstein-projective $\m$-modules.

\begin{thm} \label{mainthm} \ Let $X = (X_i, \ X_\alpha, \ i\in
Q_0, \ \alpha\in Q_1)$ be an arbitrary $\m$-module. Then
$X$ is Gorenstein-projective if and only if $X$ is a monic $\m$-module
satisfying condition ${\rm(G)}$, where

\vskip5pt

{\rm(G)} \ \ For each $i \in Q_0$, \ $X_i$ and
$X_i/ (\sum\limits_{\alpha\in \ha(\to i)}\Ima
X_\alpha)$ are Gorenstein-projective $A$-modules.
\end{thm}

Note that for a monic $\m$-module $X$, by ${\rm(m1)}$ we have $X_i/
(\sum\limits_{\alpha\in \ha(\to i)}\Ima X_\alpha) = X_i/
(\bigoplus\limits_{\alpha\in \ha(\to i)}\Ima X_\alpha).$

\subsection{}  Before a proof we give some applications.  If  $G$ is an indecomposable Gorenstein-projective $A$-module, then
$G\otimes P(v)$ is an indecomposable monic $\m$-module for each
$v\in Q_0$ (cf. Lemma \ref{adj} and Example \ref{exmmonic}$(3)$). It
is easy to see that $G\otimes P(v)$ satisfies condition {\rm(G)}.
Thus by Theorem \ref{mainthm} we have

\begin{cor} \label{GA} \ Let $G$ be an indecomposable Gorenstein-projective $A$-module. Then $G\otimes P(v)$ is an indecomposable Gorenstein-projective $\m$-module for each
$v\in Q_0$.
\end{cor}

In Example \ref{exmmonic}$(3)$ we have known that projective
$\m$-modules are monic.

\begin{cor} \label{A=k} \ Monic $\m$-modules are exactly projective $\m$-modules if and only if $A$ is a semisimple algebra.

\vskip5pt

In particular, monic $kQ/I$-modules are exactly projective
$kQ/I$-modules.
\end{cor}

\noindent{\bf Proof.} \ It suffices to prove that $A$ is semisimple
if and only if any monic $\m$-module is projective. Without loss of
generality we may assume that $A$ is {\it connected} (i.e., $A$ can
not be a product of two non-zero algebras).

\vskip5pt

Suppose that $A$ is semisimple. By Wedderburn Theorem $A\cong
M_n(R)$, the matrix algebra over a division $k$-algebra $R$, hence
$A\otimes \ kQ/I\cong M_n(R\otimes \ kQ/I)$. Since $M_n(R\otimes \
kQ/I)$ is Morita equivalent to $R\otimes \ kQ/I\cong RQ/I$,  and
${\rm gl.dim.}RQ/I < \infty$ (as in the case of $R = k$), we have
${\rm gl.dim.}(A\otimes \ kQ/I) < \infty$. Now let $M$ be a monic
$\m$-module. Then by Theorem \ref{mainthm} $M$ is a
Gorenstein-projective $\m$-module (since in  this case any
$A$-module is projective, and hence condition ${\rm(G)}$ holds
automatically). While by  [14, 10.2.3] Gorenstein-projective modules
over an algebra of finite global dimension must be projective, it
follows that $M$ is a projective $\m$-module.

\vskip5pt

Conversely, assume that any monic $\m$-module is projective. Let $M$
be an arbitrary $A$-module. Consider $\m$-module $X = M\otimes
P(1)$, where $P(1)$ is the simple projective $kQ/I$-module at sink
vertex $1$. Then $X$ is a monic $\m$-module, and hence a projective
$\m$-module by assumption. Thus $1$-st branch $M$ of $X$ is a
projective $A$-module (cf. Lemma \ref{adj}$(2)$). This proves that
$A$ is semisimple. $\s$

\vskip5pt

For a Frobenius category we refer to [27] and [21, Appendix A].

\begin{cor} \label{self} \ The following are equivalent$:$

\vskip5pt

${\rm(i)}$ \  $A$ is a self-injective algebra$;$

\vskip5pt

${\rm(ii)}$ \ $\mathcal {GP}(A\otimes kQ/I) = {\rm mon}(Q, I, A);$
\vskip5pt

${\rm(iii)}$ \ ${\rm mon}(Q, I, A)$ is a {\rm Frobenius} category.
\end{cor}
\noindent{\bf Proof.} \ ${\rm(i)} \Longrightarrow {\rm(ii)}$: If $A$
is self-injective, then every $A$-module is Gorenstein-projective,
and hence ${\rm(ii)}$ follows from Theorem \ref{mainthm}. The
implication ${\rm(ii)} \Longrightarrow {\rm(iii)}$ is well-known.

\vskip5pt

${\rm(iii)} \Longrightarrow {\rm(i)}$: Taking sink vertex $1$ of
$Q$, by Proposition \ref{injectives} $D(A_A)\otimes P(1)$ is an
injective object of ${\rm mon}(Q, I, A)$, where $D: = \Hom_k(-, k)$,
hence by assumption it is a projective object of ${\rm mon}(Q, I,
A)$. While each branch of a projective object of ${\rm mon}(Q, I,
A)$ is a projective $A$-module (cf. Theorem \ref{resolving}), it
follows that the $1$-st branch $D(A_A)$ of $D(A_A)\otimes P(1)$ is a
projective $A$-module, i.e., $A$ is self-injective. $\s$

\vskip5pt

Let $D^b(\m)$ be the bounded derived category of $\m$, and
$K^b(\mathcal P(\m))$ the bounded homotopy category of $\mathcal
P(\m)$. By definition the singularity category $D^b_{sg}(\m)$ of
$\m$ is the Verdier quotient $D^b(\m)/K^b(\mathcal P(\m))$. If $\m$
is Gorenstein, then there is a triangle-equivalence
$D^b_{sg}(\m)\cong
 \underline {\mathcal
{GP}(\m)}$  ([8, 4.4.1]; [16, 4.6]). Note that if $A$ is Gorenstein,
then $\m = A\otimes_k \ kQ/I$ is also Gorenstein. By Corollary
\ref{self} we have

\begin{cor} \label{cor} \ Let $A$ be a self-injective algebra, $Q$,  $I$ and  $\m = A\otimes \ kQ/I$ as
usual. Then there is a triangle-equivalence \ $D^b_{sg}(\m)\cong
\underline {{\rm Mon}(Q, I, A)}.$
\end{cor}

\subsection{}  For proving Theorem \ref{mainthm} we first show that the condition $(G)$ can be simplified.

\begin{lem} \label{simplerG} \ Let $X = (X_i, \ X_\alpha, \ i\in
Q_0, \ \alpha\in Q_1)$ be a monic $\m$-module such that $X_i/
(\sum\limits_{\alpha\in \ha(\to i)}\Ima X_\alpha)$ is
Gorenstein-projective for each $i \in Q_0$. Then

\vskip5pt

$(1)$ \ For each non-zero path $p$, ${\rm Im}X_p$ is
Gorenstein-projective.

\vskip5pt

$(2)$ \ $X$ satisfies condition ${\rm(G)}$.
\end{lem}

\noindent{\bf Proof.} \  $(1)$ \ Note that ${\rm Im}X_p\cong
X_{s(p)}/{\rm Ker} X_p$.  Use induction on $l_p: =\mbox{max}\{ l(q)\
|\ q\in \hp, \ e(q)=s(p)\}$. If $l_p=0$, then $s(p)$ is a source and
$X_p$ is injective, and hence ${\rm Im}X_p \cong X_{s(p)}$; while in
this case $\sum\limits_{\alpha\in \ha(\to s(p))}\Ima X_\alpha = 0$
and hence by the assumption $X_{s(p)}$ is Gorenstein-projective.

Suppose $l_p\ge 1$. Since $\mathcal {GP}(A)$ is closed under
extensions, by the assumption and the exact sequence
$$0\longrightarrow
\frac{\bigoplus\limits_{\begin{smallmatrix}\alpha\in \ha(\to s(p))
\end{smallmatrix}}\Ima
X_{\alpha}}{\Ker X_p}\longrightarrow {\rm Im}X_p\longrightarrow
\frac{X_{s(p)}}{\bigoplus\limits_{\begin{smallmatrix}\alpha\in
\ha(\to s(p))
\end{smallmatrix}} \Ima
X_{\alpha}}\longrightarrow 0$$ it suffices to prove that the left
hand side term is Gorenstein-projective.

\vskip5pt

By Theorem \ref{m4}$(2)$ $\Ker X_p = (\bigoplus\limits_{\beta\in
B_1}\Ima X_\beta) \oplus (\bigoplus\limits_{\beta\in
B_2}X_\beta(\Ker X_{p\beta})),$ where $B_1: = \{\beta \in
\mathcal{A}(\to s(p)) \ | \ \beta = {\rm la}(q) \ \mbox{for some} \
q\in K_p, \ p \beta\in I \}$ and $B_2:  = \{\beta \in
\mathcal{A}(\to s(p)) \ | \ \beta = {\rm la}(q) \ \mbox{for some} \
q\in K_p, \ p\beta\notin I\}$. Put $\Omega: = \ha(\to s(p))\setminus
(B_1\cup B_2)$.  Then
$\bigoplus\limits_{\begin{smallmatrix}\alpha\in \ha(\to s(p))
\end{smallmatrix}} \Ima X_{\alpha} = (\bigoplus\limits_{\begin{smallmatrix}\alpha\in
\Omega
\end{smallmatrix}}\Ima X_{\alpha})\oplus (\bigoplus\limits_{\begin{smallmatrix}\beta\in
B_1\cup B_2
\end{smallmatrix}} \Ima X_\beta)$ and $\Ker X_p\subseteq \bigoplus\limits_{\begin{smallmatrix}\beta\in
B_1\cup B_2
\end{smallmatrix}} \Ima X_\beta$. It follows that
\begin{align*}\frac{\bigoplus\limits_{\begin{smallmatrix}\alpha\in \ha(\to s(p))
\end{smallmatrix}} \Ima X_{\alpha}}{\Ker X_p}&=
(\bigoplus\limits_{\begin{smallmatrix}\alpha\in \Omega
\end{smallmatrix}} \Ima X_{\alpha})\oplus\frac{\bigoplus\limits_{\begin{smallmatrix}\beta\in B_1\cup B_2
\end{smallmatrix}}\Ima
X_\beta}{\Ker X_p} = (\bigoplus\limits_{\begin{smallmatrix}\alpha\in
\Omega
\end{smallmatrix}} \Ima X_{\alpha})\oplus\frac{\bigoplus\limits_{\begin{smallmatrix}\beta\in B_2
\end{smallmatrix}}\Ima
X_\beta}{\bigoplus\limits_{\beta\in B_2}X_\beta(\Ker X_{p\beta})}\\
&  \cong (\bigoplus\limits_{\begin{smallmatrix}\alpha\in \Omega
\end{smallmatrix}} \Ima X_{\alpha})\oplus \bigoplus\limits_{\begin{smallmatrix}\beta\in B_2
\end{smallmatrix}}\frac{\Ima
X_\beta}{X_\beta(\Ker X_{p\beta})}.\end{align*} If $\alpha \in
\Omega$ then $l_{ \alpha }< l_p$, and hence by induction $\Ima
X_{\alpha}$ is Gorenstein-projective. If $\beta \in B_2$ then
\begin{align*}\frac{\Ima
X_\beta}{X_\beta(\Ker X_{p\beta})}&\cong \frac{X_{s(\beta)}/\Ker
X_{\beta}}{(\Ker X_{p\beta}+\Ker X_\beta)/\Ker X_\beta}
\\& \cong \frac{ X_{s(\beta)}}{\Ker X_{p\beta}+\Ker X_\beta}
\\ & = \frac{ X_{s(\beta)}}{\Ker X_{p\beta}}
= \frac{X_{s(p\beta)}}{\Ker X_{p\beta}}\\&\cong {\rm Im}X_{p\beta},
\end{align*}
and hence it is Gorenstein-projective by induction (since
$l_{p\beta}< l_p$). This proves $(1)$.

\vskip5pt

$(2)$ \ We need to show that each $X_i$ Gorenstein-projective. This
follows from the exact sequence $0\rightarrow
\bigoplus\limits_{\begin{smallmatrix}\alpha\in \ha(\to i)
\end{smallmatrix}} \Ima X_{\alpha}\rightarrow X_i\rightarrow
X_i/ (\sum\limits_{\alpha\in \ha(\to i)}\Ima X_\alpha)\rightarrow
0$, the assumption and $(1)$,  again using the fact that $\mathcal
{GP}(A)$ is closed under extensions. $\s$

\subsection{}  Let $A$ and $B$ be rings,  $M$ an
$A$-$B$-bimodule, and
$T:=\left(\begin{smallmatrix}
A&M\\
0&B
\end{smallmatrix}\right)$ the
triangular matrix ring. We assume that $T$ is an Artin algebra ([4],
p.72), and consider finitely generated $T$-modules. Recall that a
$T$-module can be identified with a triple
$\left(\begin{smallmatrix}
X\\
Y
\end{smallmatrix}\right)_\phi$, where $X\in A$-mod, $Y\in B$-mod, and $\phi:
M\otimes_{B} Y\rightarrow X$ is an $A$-map. A $\m$-map
$\left(\begin{smallmatrix}
X \\
   Y
 \end{smallmatrix}\right)_{\phi}\rightarrow \left(\begin{smallmatrix}
   X'  \\
   Y'
 \end{smallmatrix}\right)_{\phi'}$
 can be identified with a pair
 $\left(\begin{smallmatrix}
   f  \\
   g
 \end{smallmatrix}\right)$,
where $f\in {\Hom}_A (X, X'),$ \ $g\in \Hom_B (Y, Y')$, such that
$f\phi = \phi'({\rm id}\otimes g)$. The following description of
$\mathcal {GP}(T)$ is a special case of [49, Thm. 1.4].

\begin{lem} \label{compatible}  Let $M$ be an $A$-$B$-bimodule with ${\rm proj.dim}_AM<\infty$ and ${\rm proj.dim}M_B<\infty$, and $T: = \left(\begin{smallmatrix}
A&M\\
0&B
\end{smallmatrix}\right)$.  Then
$\left(\begin{smallmatrix}X\\Y\end{smallmatrix}\right)_\phi\in
\mathcal {GP}(T)$ if and only if $\phi: M\otimes_B Y\rightarrow X$
is an injective $A$-map, $\cok\phi\in \mathcal {GP}(A)$, and $Y\in
\mathcal {GP}(B)$.
\end{lem}

In order to apply Lemma \ref{compatible} we put
\begin{align*}& \rho': = \{\rho_i\in \rho \ | \ s(\rho_i)\ne n\}; \ \ \ \ \ \
\
I': = \langle \rho_i \ | \ \rho_i\in \rho'\rangle;\\
&Q': = \mbox{the quiver obtained from} \ Q \ \mbox{by deleting
vertex} \  n;\\
& P(n): = \mbox{the indecomposable projective} \ kQ/I\mbox{-module
at} \ n;\\&\Lambda': = A\otimes kQ'/I'; \ \ \ \ \ \ \ \ \ \ \ \ \ \
\ M: = A\otimes{\rm rad} P(n).\end{align*}

\noindent Then $M$ is a $\m'$-$A$-bimodule. The point is that $\m =
A\otimes kQ/I$ is of the form $\Lambda =\left(\begin{smallmatrix}
\Lambda' &M\\
0&A
\end{smallmatrix}\right)$.
Since ${\rm gl.dim.}kQ'/I' < \infty$, we have ${\rm
proj.dim}_{kQ'/I'}{\rm rad}P(n) < \infty,$  and hence  ${\rm
proj.dim}_{\m'} M< \infty$. Clearly, as a right $A$-module $M$ is
projective. For applying Lemma \ref{compatible}, we write a
$\Lambda$-module $X= (X_i, \ X_\alpha, \ i\in Q_0, \ \alpha\in Q_1)$
as $X = \left(\begin{smallmatrix}
X'\\
X_n
\end{smallmatrix}\right)_{\phi}$,
where $X'= (X_i, \ X_\alpha, \ i\in Q'_0, \ \alpha\in Q'_1)$ is a
$\Lambda'$-module, $X_n$ is an $A$-module, and  $\phi: M\otimes_A
X_n\rightarrow X'$ is a $\Lambda'$-map. The explicit expression of
$\phi$ is given in the proof of Lemma \ref{mono1}.

\vskip5pt

By a direct translation from Lemma \ref{compatible} in this
setting-up, we have

\begin{lem} \label{mainlem} \ Let $X=
\left(\begin{smallmatrix}
X'\\
X_n
\end{smallmatrix}\right)_{\phi}$ be a $\Lambda$-module.
Then $X\in\mathcal {GP}(\m)$  if and only if $X$ satisfies the
conditions$:$

\vskip5pt

${\rm(i)}$ \ $X_n\in\mathcal {GP}(A);$

\vskip5pt

${\rm(ii)}$ \ $\phi: M\otimes_A X_n \longrightarrow X'$ is an
injective $\m'$-map$;$

\vskip5pt

${\rm(iii)}$ \ $\cok\phi\in\mathcal {GP}(\m')$.
\end{lem}

Theorem \ref{mainthm} will be proved by applying Lemma \ref{mainlem}
and using induction on $|Q_0|$. Thus, in the rest part of this paper
we write a $\Lambda$-module $X= (X_i, \ X_\alpha, \ i\in Q_0, \
\alpha\in Q_1)$ as $X = \left(\begin{smallmatrix}
X'\\
X_n
\end{smallmatrix}\right)_{\phi}$, and keep all the notations $I'$, $Q'$, $\Lambda'$, $P(n), \ M$, $X$, $X'$, $X_n$ and $\phi$.
For an integer $m\ge 0$ and a module $N$, let $N^m$ denote the
direct sum of $m$ copies of $N$ (if $m = 0$ then $N^m: = 0$).

\begin{lem} \label{mono1} \ Let $X= \left(\begin{smallmatrix}
X'\\
X_n
\end{smallmatrix}\right)_{\phi}$ be a $\Lambda$-module.
Then

\vskip5pt

$(1)$ \ $\phi: M\otimes_A X_n \rightarrow X'$ is an injective
$\m'$-map if and only if $\sum\limits_{p \in \mathcal{P}(n \to
i)}\Ima X_p =\bigoplus\limits_{p \in \hp(n \to i)}\Ima X_p$, and
$X_p$ is injective for all $i\in Q_0'$  and $p\in \mathcal P(n \to
i)$.

\vskip5pt

$(2)$ \ If  $X$ is  monic then $\phi: M\otimes_A X_n\rightarrow X'$
is injective$.$

\vskip5pt

$(3)$ \  $\cok \phi=(X_i/(\bigoplus\limits_{p\in \hp(n\rightarrow
i)}\Ima X_p), \ \widetilde{X_\alpha}, \ i\in Q'_0, \ \alpha\in
Q'_1)$, where for each $\alpha: j\rightarrow i$ in $Q'_1$,
$$\widetilde{X_\alpha}: X_j/(\bigoplus\limits_{q\in \hp(n\rightarrow
j)}\Ima X_q) \rightarrow {X_i}/ ({\bigoplus\limits_{p\in
\hp(n\rightarrow i)}\Ima X_p})$$ is the $A$-map induced by
$X_{\alpha}$. Explicitly,  \ $\widetilde{X_\alpha}(\overline {x_j})
= X_\alpha(x_j) + \bigoplus\limits_{p\in \hp(n\rightarrow i)}\Ima
X_p$, where \ $\overline {x_j}: = x_j + \bigoplus\limits_{q\in
\hp(n\rightarrow j)}\Ima X_q\in X_j/\bigoplus\limits_{q\in
\hp(n\rightarrow j)}\Ima X_q$.
\end{lem}
\noindent {\bf Proof.} $(1)$ \ For $i \in Q'_0$, put $m_i : =
|\mathcal{P}(n \to i)|$. As a representation of $(Q', I')$ over $k$,
${\rm rad} P(n)$ can be written as $\left(\begin{smallmatrix}
k^{m_1}\\
\vdots\\
k^{m_{n-1}}\\
\end{smallmatrix}\right)$, so we have isomorphisms of $\m'$-modules
$$M\otimes_A X_n \cong ({\rm rad}
P(n)\otimes_kA)\otimes_A X_n \cong {\rm rad} P(n)\otimes_k X_n \cong \left(\begin{smallmatrix} X_n^{m_1}\\
\vdots \\
X_n^{m_{n-1}}\\
\end{smallmatrix}\right).$$ Let $\hp(n\to i) = \{p_1, \cdots,p_{m_i}\}$.
Then $\phi: M\otimes_A X_n\rightarrow X'$ reads
$$\phi = \left(\begin{smallmatrix}
\phi_1\\
\vdots\\
\phi_{n-1}\\
\end{smallmatrix}\right):  \left(\begin{smallmatrix}
X_n^{m_1}\\
\vdots\\
X_n^{m_{n-1}}\\
\end{smallmatrix}\right)
\rightarrow  \left(\begin{smallmatrix}
X_1\\
\vdots\\
X_{n-1}\\
\end{smallmatrix}\right)
,$$ where $\phi_i=(X_{p_1}, \ \cdots,  \ X_{p_{m_i}}): X_n^{m_i}
\rightarrow X_i$. So $\phi$ is injective if and only if $\phi_i$ is
injective for each $i\in Q'_0$; and if and only if $\sum\limits_{p
\in \mathcal{P}(n \to i)}\Ima X_p =\bigoplus\limits_{p \in \hp(n \to
i)}\Ima X_p $, and $X_p$ is injective for all $p\in \mathcal P(n \to
i)$.

\vskip5pt

$(2)$ \ By Theorem \ref{m4}$(3)$ $\sum\limits_{p \in \mathcal{P}(n
\to i)}\Ima X_p =\bigoplus\limits_{p \in \hp(n \to i)}\Ima X_p$ for
all $i\in Q_0'$, and by Theorem \ref{m4}$(1)$ $X_p$ is injective for
all $p\in \mathcal P(n \to i)$ and $i\in Q_0'$, and hence $\phi$ is
injective by $(1)$.

\vskip5pt

$(3)$ follows from the proof of $(1)$. $\s$

\subsection{} For $i, j \in Q'_0$, we put $\mathcal P': = \{
\mbox{path of} \ Q'\}$, \   $\ha'(\to i): = \{\alpha\in Q'_1 \ | \
e(\alpha) =i\}$, \  and $\hp'(\to i): = \{\ p\in\mathcal P' \ |  \
e(p) = i,  \ l(p) \ge 1,  \
 p \notin I'\}.$

\begin{lem} \label{cokermonicandG} \ Let
$X=\left(\begin{smallmatrix}
X'\\
X_n
\end{smallmatrix}\right)_{\phi}$ be a monic $\Lambda$-module satisfying ${\rm(G)}$. Then $\cok\phi$ is a monic $\Lambda'$-module satisfying ${\rm(G)}$.
\end{lem}
\noindent {\bf Proof.} Claim 1:  \  $\cok\phi$ satisfies ${\rm
(m1)}$, i.e.,
${\sum\limits_{\alpha\in \ha'(\to i)}\Ima
\widetilde{X_{\alpha}}}={\bigoplus\limits_{\alpha\in \ha'(\to
i)}\Ima \widetilde{X_{\alpha}}}$ for each $i\in Q_0'$.

\vskip5pt

In fact, let ${\sum\limits_{\alpha\in \ha'(\to i)}
\widetilde{X_{\alpha}}}(\overline{x_{\alpha}})=0$ with each
$x_{\alpha}\in X_{s(\alpha)}$. Then
\begin{align*}  \sum\limits_{\alpha\in \ha'(\to
i)} X_\alpha (x_{\alpha}) \in &\bigoplus\limits_{p\in \hp(n\to i)}
\Ima X_p  = \sum\limits_{\beta\in \ha(n\rightarrow i)} \Ima X_\beta
+ \sum\limits_{p \in \hp (n\rightarrow i)\setminus \ha(n\to i)} \Ima
X_p \\ & =\sum\limits_{\beta \in \ha(n\rightarrow i)} \Ima X_\beta +
\sum\limits_{\alpha \in \ha(s(\alpha)\to i)} X_\alpha(\sum
\limits_{\begin{smallmatrix}q\in \hp(n\rightarrow s(\alpha))\\
\alpha q\notin I \end{smallmatrix}}\Ima X_q).\end{align*} Since
$s(\alpha) \ne n$, by ${\rm (m1)}$ on $X$ we have
$$X_\alpha(x_\alpha) \in X_\alpha(\sum \limits_{\begin{smallmatrix}q\in \hp(n\rightarrow s(\alpha))\\
\alpha q\notin I \end{smallmatrix}}\Ima X_q) \subseteq
\bigoplus\limits_{p\in \hp(n\to i)} \Ima X_p$$ for each $\alpha\in
\ha'(\to i)$. This means $\widetilde{X_\alpha}(\overline{x_\alpha})=
0.$

\vskip5pt

Claim 2:  \ $\cok\phi$ satisfies ${\rm (m2)}$, i.e.,
$\Ker\widetilde{X_{\alpha}}= {\sum\limits_{\begin
{smallmatrix}q\in \mathcal{P'}(\to s(\alpha))\\
\alpha q\in I'\end{smallmatrix}}\Ima \widetilde{X_q}}$ for each
$\alpha\in Q_1'$.

\vskip5pt

In fact, clearly $\sum\limits_{\begin {smallmatrix}q\in \mathcal{P'}(\to s(\alpha)) \\
\alpha q\in I'\end{smallmatrix}}\Ima \widetilde{X_q}\subseteq
\Ker\widetilde{X_\alpha}.$ \  Let
$\widetilde{X_{\alpha}}(\overline{x})=0$ with $x\in X_{s(\alpha)}$.
Then $X_{\alpha}(x)\in {\bigoplus\limits_{p\in \hp(n\rightarrow
e(\alpha))}\Ima X_p}.$ By the similar argument as in Claim 1
we have $X_{\alpha}(x) \in X_\alpha(\sum \limits_{\begin{smallmatrix}q\in \hp(n\rightarrow s(\alpha))\\
\alpha q\notin I \end{smallmatrix}}\Ima X_q)$. By ${\rm (m2)}$ on $X$
we have in the quotient \begin{align*}\overline x \in {\rm Ker} X_\alpha + \bigoplus\limits_{q\in \hp(n\rightarrow s(\alpha))}\Ima X_q &= \sum\limits_{\begin {smallmatrix} q'\in\hp(\to s(\alpha))\\
\alpha q'\in I\end{smallmatrix}} \Ima X_{q'}+
\bigoplus\limits_{\begin{smallmatrix}q\in \hp(n\rightarrow
s(\alpha))\end{smallmatrix}}\Ima X_q\\& = \sum\limits_{\begin {smallmatrix} q'\in\hp'(\to s(\alpha))\\
\alpha q'\in I'\end{smallmatrix}} \Ima X_{q'}+
\bigoplus\limits_{q\in \hp(n\rightarrow s(\alpha))}\Ima
X_q.\end{align*} It follows
that $\overline{x}\in\sum\limits_{\begin {smallmatrix}q'\in \mathcal{P'}(\to s(\alpha)) \\
\alpha q'\in I'\end{smallmatrix}}\Ima \widetilde{X_{q'}}$. This
proves Claim 2.

\vskip5pt

Claim 3:  \ $\cok\phi$ satisfies ${\rm (G)}$.

\vskip5pt In fact, by Lemma \ref{simplerG}$(2)$ it suffices to prove
$\frac{X_i/ (\bigoplus\limits_{p\in \hp(n\rightarrow i)}\Ima
X_p)}{\bigoplus \limits_{\alpha \in \ha'(\to i)}\Ima \widetilde
{X_\alpha}}$ is Gorenstein-projective for each $i \in Q'_0$.  Since
${\bigoplus\limits_{p \in \hp(n \to i)\setminus \ha(n \to i)}\Ima
X_p} \subseteq  \sum\limits_{\alpha \in \ha(\to i)}\Ima X_{\alpha}$,
it follows that
\begin{align*}\sum\limits_{\alpha \in \ha'(\to i)}\Ima \widetilde{X_{\alpha}} & =
\frac{{\sum\limits_{\alpha \in \ha'(\to i)}\Ima X_{\alpha} +
\bigoplus\limits_{p \in \hp(n \to
i)}\Ima X_p}}{{\bigoplus\limits_{p \in \hp(n \to i)}\Ima X_p}}\\
& = \frac{{\sum\limits_{\alpha \in \ha(\to i)}\Ima X_{\alpha} +
\bigoplus\limits_{p \in \hp(n \to i)\setminus \ha(n \to i)}\Ima
X_p}}{{\bigoplus\limits_{p \in \hp(n \to i)}\Ima X_p}}\\ & =
(\sum\limits_{\alpha \in \ha(\to i)}\Ima
X_\alpha)/(\bigoplus\limits_{p \in \hp(n \to i)}\Ima X_p) \\ &
\stackrel{{\rm (m1)}}= (\bigoplus\limits_{\alpha \in \ha(\to i)}\Ima
X_{\alpha})/(\bigoplus\limits_{p \in \hp(n \to i)}\Ima X_p).
\end{align*} Hence $$\frac{X_i/
(\bigoplus\limits_{p\in \hp(n\rightarrow i)}\Ima X_p)}{\bigoplus
\limits_{\alpha\in \ha'(\to i)}\Ima \widetilde {X_\alpha}} \cong
X_i/(\bigoplus\limits_{\alpha \in \ha(\to i)}\Ima X_{\alpha})
\eqno(4.1)$$ which is Gorenstein-projective by ${\rm(G)}$ on $X$.
$\s$

\subsection{\bf Proof of Theorem \ref{mainthm}}

Use induction on $n=| Q_0|$.

We first prove the sufficiency. Assume that $X = (X_i, \ X_\alpha, \
i\in Q_0, \ \alpha\in Q_1)=\left(\begin{smallmatrix}
X'\\
X_n
\end{smallmatrix}\right)_{\phi}$ is a monic $\m$-module satisfying
${\rm(G)}$. We need to prove that $X$ is Gorenstein-projective. The
assertion clearly holds for $n=1$. Suppose that the assertion holds
for $n-1$ \ ($n\geq 2$).  It suffices to prove that $X$ satisfies
the conditions ${\rm(i), \ (ii)}$ and ${\rm(iii)}$ in Lemma
\ref{mainlem}. In fact, the condition ${\rm(i)}$ is contained in the
condition ${\rm(G)}$; and the condition ${\rm(ii)}$ follows from
Lemma \ref{mono1}$(2)$. By Lemma \ref{cokermonicandG} $\cok\phi$ is
a monic $\Lambda'$-module satisfying ${\rm(G)}$. Since $|Q_0'| =
n-1$, it follows from induction that $\cok\phi$ is
Gorenstein-projective, i.e., the condition ${\rm(iii)}$ in Lemma
\ref{mainlem} is also satisfied. This proves that $X$ is
Gorenstein-projective.

\vskip5pt

Now we prove the necessity. Assume that $X = (X_i, \ X_\alpha, \
i\in Q_0, \ \alpha\in Q_1)=\left(\begin{smallmatrix}
X'\\
X_n
\end{smallmatrix}\right)_{\phi}$ is a Gorenstein-projective $\Lambda$-module.
We need to prove that $X$ is a monic $\Lambda$-module satisfying
${\rm(G)}$. The assertion clearly holds for $n=1$. Suppose that the
assertion holds for $n-1$ \ ($n\geq 2$). By Lemma \ref{mainlem}
$X_n$ is a Gorenstein-projective $A$-module, \ $\phi: M\otimes_A X_n
\rightarrow X'$ is an injective $\m'$-map, and $\cok\phi$ is a
Gorenstein-projective $\m'$-module. Since $\phi$ is injective, by
Lemma \ref{mono1}$(1)$  we have

\vskip5pt

$(1)$ \ For each $i\in Q_0'$ there holds
$\sum\limits_{p\in \hp(n\to i)}\Ima  X_p=\bigoplus\limits_{\begin{smallmatrix} p\in \hp(n\to i)\end{smallmatrix}}\Ima  X_p;$ and

\vskip5pt

$(2)$ \ $X_p$ is an injective $A$-map for $p\in \mathcal P(n \to i)$
and $i\in Q_0'$.

\vskip5pt

Since $\cok\phi$ is a Gorenstein-projective $\m'$-module and $|Q'_0|
= n-1$, by induction $\cok\phi$ is a monic $\Lambda'$-module
satisfying ${\rm(G)}$. Thus we have

\vskip5pt

$(3)$ \ For each $i\in Q_0'$ there holds $\sum\limits_{\alpha\in
\ha'(\to i)}\Ima \widetilde{X_\alpha} =\bigoplus\limits_{\alpha\in
\ha'(\to i)}\Ima \widetilde{X_\alpha};$ and

\vskip5pt

$(4)$ \ For each $\alpha\in Q_1'$ there holds
$\Ker\widetilde{X_{\alpha}}= \sum\limits_{\begin
{smallmatrix}q\in \mathcal{P'}(\to s(\alpha))\\
\alpha q\in I'\end{smallmatrix}}\Ima \widetilde{X_q}.$

\vskip5pt

{\bf Claim 1}:  \ $X$ satisfies ${\rm (m1)}$, i.e., for each $i\in
Q_0$ there holds $\sum\limits_{\begin {smallmatrix} \alpha\in
\ha(\to i)
\end{smallmatrix}}\Ima
X_\alpha = \bigoplus\limits_{\begin {smallmatrix} \alpha\in \ha(\to
i) \end{smallmatrix}}\Ima X_\alpha.$

\vskip5pt
In fact, suppose that
$\sum\limits_{\begin {smallmatrix} \alpha\in \ha(\to i)
\end{smallmatrix}}X_{\alpha}(x_\alpha) =0$
with each $x_\alpha\in X_{s(\alpha)}$. Then $\sum\limits_{\alpha\in
\ha'(\to i)}\widetilde{X_{\alpha}}(\overline{x_\alpha}) =0.$ By
$(3)$ we have $\widetilde{X_{\alpha}}(\overline{x_\alpha})=0$ for
each $\alpha\in \ha'(\to i)$, and hence  $\overline{x_\alpha}\in
\sum\limits_{\begin
{smallmatrix}q\in \mathcal{P'}(\to s(\alpha))\\
\alpha q\in I'\end{smallmatrix}}\Ima \widetilde{X_q}$  by $(4)$.  So
there exists $y_\alpha\in \sum\limits_{\begin
{smallmatrix}q\in \mathcal{P'}(\to s(\alpha))\\
\alpha q\in I'\end{smallmatrix}}\Ima X_q$ such that $x_\alpha -
y_\alpha\in \bigoplus\limits_{\begin{smallmatrix} p\in \hp(n\to
s(\alpha))\end{smallmatrix}}\Ima  X_p.$ Thus there exists $x'_\alpha
= \sum\limits_{\begin{smallmatrix} p\in \hp(n\to
s(\alpha))\end{smallmatrix}}X_p(x_{p, \alpha})$ with each $x_{p,
\alpha}\in X_{s(p)}$ such that $x_\alpha - y_\alpha = x_\alpha'$.
Since $X_\alpha X_q = 0$ for $\alpha q\in I'$, we have
$X_\alpha(y_\alpha) = 0$, and hence $X_{\alpha}(x_\alpha)=
X_{\alpha}(x'_\alpha).$ It follows that
\begin{align*}0&=\sum\limits_{\begin {smallmatrix} \alpha\in \ha(\to i)
\end{smallmatrix}}X_{\alpha}(x_\alpha) = \sum\limits_{\alpha\in \ha(n \to i)}X_{\alpha}(x_\alpha) + \sum\limits_{\alpha\in \ha'(\to i)}X_{\alpha}(x_\alpha)
\\& =\sum\limits_{{\alpha}\in \ha(n \to i)}X_{\alpha}(x_\alpha)+ \sum\limits_{\alpha\in \ha'(\to i)}X_{\alpha}(\sum\limits_{p\in \hp(n\to s(\alpha))}X_p(x_{p, \alpha})).\end{align*}
By $(1)$ this sum is a direct sum, hence
 $X_{\alpha}(x_\alpha)= 0$ for  $\alpha\in \ha(n \to i)$, and \ $X_{\alpha}X_p(x_{p, \alpha}) = 0$ for
 $\alpha\in \ha'(\to i)$ and $p\in \hp(n\to s(\alpha))$. Thus all together $X_{\alpha}(x_\alpha)=0$ for each $\alpha\in \ha(\to i)$. This proves Claim 1.

\vskip5pt

{\bf Claim 2}: \ $X$ satisfies ${\rm (m2)}$, i.e., for each
$\alpha\in Q_1$ there holds $\Ker X_\alpha
=\sum\limits_{\begin {smallmatrix} q\in\hp(\to s(\alpha))\\
\alpha q\in I\end{smallmatrix}} \Ima X_{q}.$

\vskip5pt

Let $x\in \Ker X_\alpha.$ Then $\overline{x}\in \Ker
\widetilde{X_{\alpha}}.$ By $(4)$ we have
$$\overline{x}=\sum\limits_{\begin
{smallmatrix}q\in \mathcal{P'}(\to s(\alpha))\\
\alpha q\in I'\end{smallmatrix}}
   X_q(x_q)+ \bigoplus\limits_{\begin {smallmatrix}p\in \hp(n\to s(\alpha))\end{smallmatrix}}
 \Ima X_p$$
for some $x_q\in X_{s(q)}$,   and hence there is $y_p\in X_n$  for
each $p\in\hp(n\to s(\alpha))$, such that
\begin{align*} x &=\sum\limits_{\begin
{smallmatrix}q\in \mathcal{P'}(\to s(\alpha))\\
\alpha q\in I'\end{smallmatrix}}
 X_q(x_q)+ \sum\limits_{\begin {smallmatrix}p\in \hp(n\to s(\alpha)) \end{smallmatrix}} X_p(y_p)
 \\ & =\sum\limits_{\begin
{smallmatrix}q\in \mathcal{P'}(\to s(\alpha))\\
\alpha q\in I'\end{smallmatrix}}
 X_q(x_q)+ \sum\limits_{\begin {smallmatrix}p\in \hp(n\to s(\alpha))\\ \alpha p\in I\end{smallmatrix}} X_p(y_p)
 + \sum\limits_{\begin {smallmatrix}p\in \hp(n\to s(\alpha))\\ \alpha p\notin I\end{smallmatrix}} X_p(y_p)\end{align*}
Since $X_\alpha(x)=0$ and $X_\alpha X_q=0$ for $\alpha q\in I$, we
have $\sum\limits_{\begin {smallmatrix}p\in \hp(n\to s(\alpha))\\
\alpha p\notin I\end{smallmatrix}} X_\alpha X_p(y_p)= 0.$ By $(1)$
we have $X_{\alpha}X_p(y_p)=0$ for each $p\in \hp(n\to s(\alpha))$
with $\alpha p\notin I$. But by $(2)$ $X_\alpha X_p$ is injective,
so $y_p = 0$ for all $p\in \hp(n\to s(\alpha))$ with $\alpha p\notin
I$. Thus $$x=\sum\limits_{\begin
{smallmatrix}q\in \mathcal{P'}(\to s(\alpha))\\
\alpha q\in I'\end{smallmatrix}}
 X_q(x_q)+ \sum\limits_{\begin {smallmatrix}p\in \hp(n\to s(\alpha))\\ \alpha p\in I\end{smallmatrix}} X_p(y_p) \in
\sum\limits_{\begin {smallmatrix} q\in\hp(\to s(\alpha))\\
\alpha q\in I\end{smallmatrix}} \Ima X_{q}.$$ This proves Claim 2.

\vskip5pt

By {\bf Claims 1} and {\bf 2} $X$ is a monic $\m$-module. It remains
to prove that  $X$ satisfies ${\rm(G)}$. By Lemma
\ref{simplerG}$(2)$ it suffices to prove that  $X_i/
(\bigoplus\limits_{\alpha\in \ha(\to i)}\Ima X_\alpha)$ is a
Gorenstein-projective $A$-module for each $i\in Q_0$. We have
claimed that $\cok \phi$ satisfies ${\rm(G)}$, i.e., for each $i\in
Q_0'$, \ $(X_i/(\bigoplus\limits_{p \in \hp(n \to i)}\Ima
X_p))/(\bigoplus\limits_{\alpha\in \ha(\to i)} \Ima \widetilde
{X_\alpha})$ is a Gorenstein-projective $A$-module, which is
isomorphic to $X_i/ (\bigoplus\limits_{\alpha\in \ha(\to i)}\Ima
X_\alpha)$ by $(4.1)$. This completes the proof of Theorem
\ref{mainthm}. $\s$

\vskip20pt

Xiu-Hua Luo, \ xiuhualuo2014$\symbol{64}$163.com

Dept. of Math., \ Nantong Univ., \  Nantong 226019,
 \  China

\vskip5pt

Pu Zhang, \ pzhang$\symbol{64}$sjtu.edu.cn

Dept. of Math., \ Shanghai Jiao Tong Univ., \ Shanghai 200240, \
China

\end{document}